\documentclass[11pt]{article}
\usepackage{epsfig}
\usepackage{latexsym,graphicx}
\title{\bf A finite dimensional approach to Bramham's approximation theorem}

\author{ Patrice LE CALVEZ \thanks{ 
Institut de Math\'ematiques de Jussieu, UMR 7586 CNRS. Universit\'e Pierre et Marie Curie. Case 247, 4 Place Jussieu,
75005 Paris Cedex, France. E-mail: lecalvez@math.jussieu.fr
} }

 \date{July 19 th, 2013}

\def\picture #1 by #2 (#3){
  \vbox to #2{
    \hrule width #1 height 0pt depth 0pt
    \vfill
    \special{picture #3} 
    }
  }

\def\scaledpicture #1 by #2 (#3 scaled #4){{
  \dimen0=#1 \dimen1=#2
  \divide\dimen0 by 1000 \multiply\dimen0 by #4
  \divide\dimen1 by 1000 \multiply\dimen1 by #4
  \picture \dimen0 by \dimen1 (#3 scaled #4)}
  }

\def\R{{\bf R}}
\def\Z{{\bf Z}}

\def\Q{{\bf Q}}
\def\T{{\bf T}}

\def\D{{\bf D}}
\def\S{{\bf S}}

\def\e{\varepsilon}
\begin{document}
\def\eqalign#1{\null\,\vcenter{\openup\jot
\ialign{\strut\hfil$\displaystyle{##}$&
$\displaystyle{{}##}$\hfil \crcr #1\crcr }}\,}

\def\eqalignno#1{\displ@y \tabskip=\@centering
\halign to\displaywidth{\hfil$\@lign\displaystyle{##}$
\tabskip=0pt &$\@lign\displaystyle{{}##}$

\hfil\tabskip=\@centering
$\llap{$\@lign##$}\tabskip=Opt\crcr #1\crcr}}

\thispagestyle{empty}
\maketitle

\noindent{\bf Abstract:} \enskip Using pseudoholomorphic curves techniques from symplectic geometry, Barney Bramham proved recently that every smooth irrational  
 pseudo-rotation of the unit disk is the limit, for the $C^0$ topology, of a sequence of smooth periodic diffeomorphisms. We give here a finite dimensional proof of this result that works in the case where the pseudo-rotation is smoothly conjugate to a rotation on the boundary circle. The proof extends to $C^1$ pseudo rotations and is based on the dynamical study of the gradient flow associated to a generating family of functions given by Chaperon's broken geodesics method.



\bigskip
\bigskip
\bigskip
\bigskip

\noindent{\large {\bf 0. Introduction}}

\bigskip

We will denote by $\D$ the closed unit disk of the euclidean plane and by $\S$ the unit circle. An {\it irrational pseudo-rotation}  is an area preserving homeomorphism of $\D$ that fixes $0$ and that does not possess any other periodic point. To such a homeomorphism is associated an irrational number $\overline{\alpha}\not\in\Q/\Z$ characterized by the following: every point admits $\overline{\alpha}$ as a rotation number. To give a precise meaning to this sentence, choose a lift $\widetilde f$  of $f\vert_{\D\setminus\{0\}}$ to the universal covering space $\widetilde \D=\R\times(0,1]$. There exists $\alpha\in\R$ satisfying $\alpha+\Z=\overline{\alpha}$ such that for every compact set $\Xi\subset\D\setminus\{0\}$, and every $\varepsilon >0$, one can find $N\geq 1$ such that $$n\geq N \enskip\mathrm{and}\enskip \widetilde z\in \pi^{-1}(\Xi)\cap\widetilde f^{-n}(\pi^{-1}(\Xi))\Rightarrow \left\vert {p_1(\widetilde f^n(\widetilde z))-p_1(\widetilde z)\over n}-\alpha\right\vert \leq \varepsilon,$$ where $\pi:  (r,\theta)\mapsto (r\cos2\pi\theta, r\sin 2\pi\theta)$ is the covering projection and
$p_1: (r,\theta)\mapsto r$ the projection on the first factor. In particular the Poincar\'e rotation number of $f\vert_{\S}$ is $\overline{\alpha}$. In the case where $f$ is a $C^k$ diffeomorphism, $1\leq k\leq \infty$, we will say that $f$ is a $C^k$ irrational pseudo-rotation.
Constructions of dynamically interesting irrational pseudo-rotations are based on the method of fast periodic approximations, starting from the seminal paper of Anosov and Katok \cite{AK} (see \cite {FatHe}, \cite{FayK}, \cite{FayS}, \cite{Ha} for further developments about this method, see \cite{BeCrLerP}, \cite{BeCrLer}  for other results on irrational pseudo-rotations).

\medskip

 Barney Bramham has recently proved the following (see \cite{Br1}):

\bigskip
\noindent{\sc Theorem 1:}\enskip\enskip {\it Every $C^{\infty}$ irrational pseudo-rotation $f$ is the  limit, for the $C^0$ topology, of a sequence of periodic $C^{\infty}$ diffeomorphisms.}

\bigskip
The result is more precise. Let $(q_n)_{n\geq 0}$ be a sequence of positive integers such that the sequence $(q_n\overline{\alpha})_{n\geq 0}$ converges to $0$ in $\T^1=\R/\Z$. One can construct a sequence  of homeomorphisms $(f_n)_{n\geq 0}$ fixing $0$ and satisfying $(f_n)^{q_n}=\mathrm{Id}$ that converges to $f$ in the $C^0$ topology. Such a map $f_n$ is $C^0$ conjugate to a rotation of rational angle (mod. $\pi$). Approximating the conjugacy by a $C^{\infty}$ diffeomorphism permits to approximate $f_n$ by a $C^{\infty}$ diffeomorphism of same period. 

\medskip

The proof uses pseudoholomorphic curves techniques from symplectic geometry. Trying to find a finite dimensional proof of this result is natural, as some results of symplectic geometry admit finite dimensional proofs by the use of generating families. A seminal example is Chaperon's proof of Conley-Zehnder's Theorem via broken geodesics method (see \cite{Ch}): if $F$ is the time one map of a Hamiltonian flow on the torus $\T^{2r}=\R^{2r}/\Z^{2r}$, a function can be constructed on a space $\T^{2r}\times \R^{2n}$ whose critical points are in bijection with the contractible fixed points of $F$. Studying the dynamics of the gradient vector field $\xi$ permits to minimize the number of critical points. Writing $F$ as a composition of diffeomorphisms $C^1$ close to the identity is the way Chaperon constructs a generating family. Decomposing $F$ in monotone twist maps alternatively positive or negative is another possible way. It is the fact that $F$ is isotopic to the identity that is essential in the construction of the vector field $\xi$, but in the general case $\xi$ has no reason to be a gradient vector field and its dynamics may be more complex. Nevertheless, if $r=1$ the vector field will satisfy some ``canonical dissipative properties'' and its dynamics can be surprisingly well understood (see \cite{Lec2} for the case where $F$ is decomposed in monotone twist maps).  Among the applications, one  can note the following approximation result (see \cite{Lec3}):  every minimal $C^1$ diffeomorphism $F$ of $\T^2$ that is isotopic to the identity is a limit in the $C^0$ topology of a sequence of periodic diffeomorphisms. The proofs given in \cite{Br1} and \cite {Lec3} share a thing in common: the construction of a foliation satisfying a  certain ``dynamical transverse property'' on which a finite group acts, the approximating map being naturally related to this action.  In  \cite{Br1} the foliation is defined on  $\R\times {\bf D}\times \T^1$ and the leaves are either pseudoholomorphic cylinders or pseudoholomorphic half cylinders transverse to the boundary; in \cite{Lec3}, the foliation is singular and naturally conjugate to the foliation by  orbits of $\xi$ on an invariant torus. Therefore it is natural to look for a proof of Bramham's theorem by a method close to the one given in \cite{Lec3}. The original proof of Theorem 1 is divided in two cases: the case where the restriction of $f$ to $\S$ is smoothly conjugate to rotation, and the case where it is not. We succeeded to treat the first case, with some improvements due to the fact that we work in the $C^1$ category but unfortunately  could not get the general case. Therefore we will prove:

\bigskip
\noindent{\sc Theorem 1':}\enskip\enskip {\it Every $C^{1}$ irrational pseudo-rotation $f$, whose restriction to $\S$ is $C^1$ conjugate to a rotation, is the limit, in the $C^0$ topology, of a sequence of periodic smooth diffeomorphisms. }

\bigskip
Observe that it is sufficient to prove Theorem 1' in the case where the restriction to $\S$ is a rotation. Indeed, every $C^1$ diffeomorphism of $\S$ can be extended to a $C^1$ area preserving diffeomorphism  of $\D$ (see \cite{BoCrW} for example). So, every $C^{1}$ irrational pseudo-rotation, whose restriction to $\S$ is $C^1$ conjugate to a rotation, is itself conjugate to a $C^{1}$ irrational pseudo-rotation, whose restriction to $\S$ is a rotation.

\bigskip
Let us explain the ideas of the proof. The first difficulty arises from the fact that $f$ is defined on a surface with boundary. If one supposes that $f\vert_{\bf S}$ is a rotation, one can extend easily our map to the whole plane. Inside a small neighborhood of ${\bf D}$ we extend our map by an integrable polar twist map and outside by a rotation whose angle is irrational (mod. $\pi$) and close but different from $2\pi\alpha$. This implies that ${\bf S}$ is accumulated by invariant circles $S_{p/q}$ on which the map is periodic with a rotation number $p/q$ that is a convergent of $\alpha$, where $\alpha+\Z=\overline{\alpha}$. Our extended map is piecewise $C^1$ and one can construct a generating family of functions that are $C^1$ with Lipschitz derivatives (see Section 1). We could have chosen to decompose $f$ in monotone twist maps in order to apply directly the results of \cite{Lec1} and \cite{Lec2},  we have preferred to use a decomposition in maps close to the identity like in \cite{Ch} to underline the fact that the way we construct the generating family is not important. One knows that for every $q\geq 1$, the fixed point set of $f^q$ corresponds to the singular points of a gradient vector field $\xi_q$ defined on a space $E_q$ depending on $q$. In particular each circle $S_{p/q}$ corresponds to a curve $\Sigma_{p/q}$ of singularities of $\xi_q$. In Section 2 we will recall the immediate properties of $\xi_q$, in particular its invariance by the natural action of $\Z/q\Z$ on $E_q$. A crucial point is the fact that $\xi_q$ is $A$ Lipschitz with a constant $A$ that does not depend on $q$. An important consequence is an uniform inequality between the $L^2$ norm of an orbit (the square root of the energy) and its $L^0$ norm.  In Section 3 we give the proofs of Theorem 1'. The fundamental result (Proposition 2) is the fact that $\Sigma_{p/q}$ bounds a disk $\Delta_{p/q}$ that contains the singular point corresponding to the fixed point $0$ and that is invariant by the flow and by the $\Z/q\Z$ action. Moreover the dynamics on $\Delta_{p/q}$ is North-South and the non trivial orbits have the same energy. This energy can be explicitely computed and is small if $p/q$ is a convergent of $\alpha$. Consequently the vector field is uniformly small on $\Delta_{p/q}$. The approximation map will be related to $\xi_q\vert_{\Delta_{p/q}}$, as it is done in \cite{Lec3}. It must be noticed that the arguments of this section are the finite dimensional analogous of the arguments of \cite{Br1}. The rest of the paper is devoted to the proof of Proposition 2. The vector field would have been $C^1$, one could have use the results of \cite{Lec2} and the following remark: the set $\{0\}\cup S_{p/q}$ is a {\it maximal unlinked fixed points set} of $f^q$, which means that there exists an isotopy from identity to $f^q$ that fixes every point of $\{0\}\cup S_{p/q}$ and there is no larger subset of the fixed point set of $f^q$ that satisfies this property. The vector field being Lipschitz, one must adapt what is known in the $C^1$ case to this wider situation. In section 4 we recall the existence a {\it dominated structure} by presenting a canonical filtration on the product flow on $E_q\times E_q$, postponing the technical proofs in the Appendix. In section 5 we explain how $\Delta_{p/q}$ appears has an ``invariant manifold'' of this dominated structure.  We have tried to write the article as much self-contained as possible.

\medskip
I would like to thank Barney Bramham for instructive and useful conversations.

\bigskip
\bigskip
\vfill\eject

\noindent{\large {\bf 1. Extension and decomposition of a pseudo-rotation}}

\bigskip
Let $f$ be an orientation preserving homeomorphism of the Euclidean plane. We will say that $f$  is {\it untwisted} if the map
$$(x,y)\mapsto(p_1(f(x,y)), y)$$ is a homeomorphism, which means that there exist two continuous functions $g, g'$ on $\R^2$ such that
$$f(x,y)=(X,Y)\Leftrightarrow\cases{ x=g(X,y)\cr Y=g'(X,y)\cr}.$$ In this case, the maps $X\mapsto g(X,y)$ and $y\mapsto g'(X,y)$ are orientation preserving homeomorphisms of $\R$. If moreover, $f$ is area preserving, the continuous form $x dy+Y dX$ is exact: there exists a  $C^1$ function $h: \R^2\to \R$ such that 
$$\displaystyle g={\partial h\over \partial y}, \enskip\displaystyle g'={\partial h\over \partial X}.$$The function $h$, defined up to an additive constant is a  {\it generating function} of $f$. 

\medskip
We will be interested in untwisted homeomorphisms satisfying some Lipschitz conditions. 
Let $f$ be an orientation preserving homeomorphism of the euclidean plane and $K\geq 1$. We will say that $f$ is a {\it $K$ Lipschitz untwisted homeomorphism } if

\smallskip
\noindent {\bf i)}\enskip \enskip $f$ is untwisted;

\smallskip
\noindent {\bf ii)}\enskip \enskip $f$ is $K$ bi-Lipschitz;

\smallskip
\noindent {\bf ii)}\enskip \enskip the maps $X\mapsto g(X,y)$ and $y\mapsto g'(X,y)$ are $K$ bi-Lipschitz;

\smallskip
\noindent {\bf iii)} \enskip \enskip the maps $y\mapsto g(X,y)$ and $x\mapsto g'(X,y)$ are $K$ Lipschitz.

\medskip

Let $f$ be a $C^1$ diffeomorphism of $\R^2$ and denote by $$\mathrm{Jac}(f)=\left(\matrix {{\partial X\over\partial x}&{\partial X\over\partial y}\cr {\partial Y\over\partial x}&{\partial Y\over\partial y}\cr}\right)$$ its Jacobian matrix. One can verify that $f$ is a $K$ Lipschitz untwisted homeomorphism  if and only if the eigenvalues of the matrix $\mathrm{Jac}(f)^t\mathrm{Jac}(f)$ lie between $K^{-2}$ and $K^2$ and if the following conditions are fulfilled

$$K^{-1}\leq{\partial X\over \partial x}\leq K, \enskip K^{-1}\leq \mathrm{det}(\mathrm{Jac}(f))^{-1}{\partial X\over\partial x}\leq K$$
and 
$$\left\vert{\partial X\over\partial y}\right\vert \leq K{\partial X\over\partial x}, 
\enskip \left\vert{\partial Y\over\partial x}\right\vert \leq K{\partial X\over\partial x}.$$
In particular all these conditions are satisfied if $\mathrm{Jac}(f)$ is (uniformly) sufficiently close to the identity matrix.

Until the end of Section 4 we suppose given a $C^1$ pseudo-rotation of rotation number $\overline\alpha$ that coincides with a rotation on $\S$. We choose a real representant $\alpha$ of $\overline\alpha$ and an irrational number $ \beta> \alpha$ such that $(\alpha, \beta)\cap\Z=\emptyset$. We extend our map to a homeomorphism $f$ of the whole plane defined in polar coordinates as follows:$$
 f(r,\theta)= \cases{(r, \theta+2\pi(\alpha+r-1))& if \enskip $r\in[1,1+\beta-\alpha]$,\cr
 (r, \theta+2\pi\beta)& if \enskip $r\geq 1+\beta-\alpha$.\cr}
$$
We get a piecewise $C^1$ area preserving transformation that satisfies the following properties:

\smallskip
\noindent-\enskip\enskip $0$ is the unique fixed point of $f$;

\smallskip
\noindent-\enskip\enskip there is no periodic point of period $q$ if $(q \alpha, q \beta)\cap\Z=\emptyset$;

\smallskip
\noindent-\enskip\enskip if $(q \alpha, q \beta)\cap\Z=\emptyset$, the set of periodic points of period $q$ can be written $\bigcup_{\alpha<p/q<\beta} S_{p/q}$, where $S_{p/q}$ is the circle of center $0$ and radius $1+p/q- \alpha$.

\bigskip
\noindent{\sc Proposition 1:}\enskip\enskip {\it For every $K>1$, one can find a decomposition $f=f_m\circ\dots \circ f_1$, where each $f_i$ is an area preserving $K$ Lipschitz untwisted homeomorphim that fixes $0$ and induces a rotation on every circle of origin ${0}$ and radius $r\geq 1$.}

\bigskip
\noindent{\it Proof.} \enskip Let denote by $f'$ the plane homeomorphism defined in polar coordinates as follows:
$$
f'(r,\theta)= \cases{(r, \theta+2\pi\alpha)& if \enskip $r\in[0,1]$,\cr
 (r, \theta+2\pi( \alpha+r-1))& if \enskip $r\in[1,1+\beta-\alpha]$,\cr
 (r, \theta+2\pi  \beta)& if \enskip $r\geq 1+\beta-\alpha$.\cr}
$$
One gets $f'=(f'_{m'}{})^{m'}$, where $$
 f'_{m'}(r,\theta)= \cases{(r, \theta+2\pi\alpha/m')& if \enskip $r\in[0,1]$,\cr
 (r, \theta+2\pi(\alpha+r-1)/m')& if \enskip $r\in[1,1+\beta-\alpha]$,\cr
 (r, \theta+2\pi  \beta/m')& if \enskip $r\geq 1+\beta-\alpha.$\cr}
$$

If $m'$ is large enough, $\widetilde f'_{m'}$ is an area preserving $K$ Lipschitz untwisted homeomorphism. Indeed it induces a diffeomorphism on each surface of equation $$r\leq 1, \enskip 1\leq r\leq 1+ \beta-\alpha, \enskip 1+\beta-\alpha\leq r$$ and the Jacobian at every point is uniformly close to the identity..

One can write $f=f''\circ f'$ where $f''$ coincides with the identity outside $\D$ and is an area preserving $C^1$ diffeomorphism, when restricted to $\D$. To get the proposition one can use the path-connectedness of the group $\mathrm{Diff}^1_{**}({\bf D})$ of area preserving $C^1$ diffeomorphisms of $\bf D$ that fix $0$ and every point of $\bf S$, when munished with the $C^1$ topology\footnote{Having not be able to find a written proof of the path-connectedness of $\mathrm{Diff}^1_{**}({\bf D})$ in the litterature, we have written one in the Appendix (Lemma 7)}.  Indeed, for every neighborhood $\cal U$  of the identity in $\mathrm{Diff}^1_{**}({\bf D})$, one can write
$f''\vert_{\D}=f''_{m''}\circ\dots \circ f''_{1}$  where $f''_i\in \cal U$. Choosing $\cal U$ sufficiently close to the identity and  extending $f''_i$ with the identity map outside $\D$, one gets a decomposition of $f''$ in area preserving $K$ Lipschitz untwisted homeomorphisms. It remains to write $m=m'+m''$ and to define
$$
 f_i= \cases{f'_{m'}& if \enskip $i\leq m'$,\cr
f''_{i-m'}& if \enskip $m'<i\leq m$.\cr}
$$
 
  \hfill$\Box$

\bigskip
\noindent{\bf Remarks.}\enskip 

\medskip
\noindent{\bf 1.}\enskip\enskip To each map $f_i$ is naturally  associated an isotopy $(f_{i,t})_{t\in[0,2]}$ starting from the identity: writing $f_i(x,y)=(X,Y)$, one sets 
 $$f_{i,t}(x,y)= \cases {((1-t) x+t X, y) & if $t\in[0,1]$,\cr
 (X, (2-t)y+ (t-1) Y) & if $t\in[1,2]$.\cr}$$ One gets an isotopy $(f_t)_{t\in[0,2m]}$ joining the identity to $f$ and fixing the origin by writing
 $$f_t= f_{1+i,t-2i}\circ f_{i}\circ \dots \circ f_1 \enskip\mathrm{if}\enskip t\in [2i, 2i+2].$$This isotopy can be lifted on the universal cover $\R\times(0,+\infty)$ of ${\bf R^2}\setminus\{0\}$ to an isotopy from the identity to a certain lift $\widetilde f$ of $f\vert_{{\bf R}^2\setminus\{0\}}$. The real rotation number defined by the restriction of $\widetilde f$ to the universal cover $\widetilde {\bf D}={\bf R}\times(0,1)$ of ${\bf D}\setminus\{0\}$ is $\alpha$.

\medskip
\noindent{\bf 2.}\enskip\enskip By choosing $m'$ sufficiently large in the proof of Proposition 1, one can suppose that for every $i\in\{1,\dots, m\}$ and for every $r\geq 1$, the rotation that coincides with $f_i$ on the circle of origin $0$ and radius $r$  is a $K$ Lipschitz untwisted homeomorphim. This fact will be used in section 4.

\bigskip
\bigskip

\noindent{\large {\bf 2. The generating family and the gradient flow}}

\bigskip
We fix $K>1$ and a decomposition $f=f_m\circ\dots \circ f_1$ given by Proposition 1. We define two families $(g_i)_{1\leq i\leq m}$, $(g'_i)_{1\leq i\leq m}$ of continuous maps as follows
$$f_i(x,y)=(X,Y)\Leftrightarrow\cases{ x=g_i(X,y),\cr Y=g'_i(X,y),\cr}$$ and a family $(h_i)_{1\leq i\leq m}$ of $C^1$ maps, such that$$\displaystyle g_i={\partial h_i\over \partial y}, \enskip\displaystyle g'_i={\partial h_i\over \partial X}.$$We extend the families 
$$(f_i)_{1\leq i\leq m}, \enskip (g_i)_{1\leq i\leq m},  \enskip (g'_i)_{1\leq i\leq m},  \enskip (h_i)_{1\leq i\leq m},$$ to $m$ periodic families $$(f_i)_{i\in\Z},  \enskip (g_i)_{i\in\Z},  \enskip (g'_i)_{i\in\Z},  \enskip (h_i)_{i\in\Z},$$

\medskip 
We fix in this section an integer $q\geq 2$ such that $(q\alpha, q\beta)\cap {\Z}\not=\emptyset$. To lighten the notations, unlike in the introduction we do not refer to $q$ while defining objects. We consider the finite dimensional vector space
$$	E= \left\lbrace {\bf z}=(z_{i})_{i\in {\Z}}\in (\R^2)^{\Z}\enskip 
\enskip\vert \enskip z_{i+mq}=z_{i},\enskip 
\mathrm{ for\enskip all}\enskip i\in {\Z}\right\rbrace
$$
munished with the scalar product
$$\left\langle
(z_{i})_{i\in {\Z}}, (z'_{i})_{i\in {\Z}}\right\rangle=\sum_{0< i\leq mq} x_ix'_i+y_iy'_i,$$ where $z_{i}=(x_i,y_i)$ and $z'_{i}=(x'_i,y'_i)$. We denote by $\Vert\enskip\Vert$ the associated Euclidean norm and write $$d({\bf z}, Z)=\inf_{{\bf z}'\in Z}  \Vert {\bf z}-{\bf z'}\Vert$$ for the distance of a point ${\bf z}$ to a set $Z\subset E$.

\medskip
We define on $E$ a vector field $\xi=(\xi_i)_{i\in\Z}$ by writing
$$
{\bf \xi}
_{i}({\bf z})=\left(y_i-g'_{i-1}(x_{i},y_{i-1}), \enskip x_i-g_i(x_{i+1},y_i)\right).$$ Note that $\xi$ is invariant by the ($q$ periodic) shift
$$\eqalign{\varphi : E&\to E\enskip ,\cr
(z_{i})_{i\in {\Z}}&\mapsto (z_{i+m})_{i\in 
{\Z}} .\cr}$$
In this section, we will state easy facts about $\xi$.
 
\bigskip
\noindent{\sc Lemma 1:}\enskip\enskip {\it The vector field is $A$ Lipschitz, where $A=\sqrt{6K^2+3}$.}

\bigskip
\noindent{\it Proof.} \enskip For every ${\bf z}=(z_i)_{i\in\Z}$ and ${\bf z}'=(z'_i)_{i\in\Z}$, one has
$$\left\vert (y_i-g'_{i-1}(x_{i},y_{i-1})) - (y'_i-g'_{i-1}(x'_{i},y'_{i-1}))\right\vert \leq \vert y_i-y'_i\vert+K\vert x_i-x'_i\vert +K \vert y_{i-1}-y'_{i-1}\vert$$ 
and
$$\left\vert (x_i-g_{i}(x_{i+1},y_{i})) - (g_{i}(x'_i-x'_{i+1},y'_{i}))\right\vert \leq \vert x_i-x'_i\vert+K\vert x_{i+1}-x'_{i+1}\vert +K \vert y_{i}-y'_{i}\vert.$$
By Cauchy-Schwarz inequality, one knows that $$(a+b +c)^2\leq 3(a^2+b^2+c^2),$$ which implies that
$$\Vert \xi_i({\bf z})-\xi_i({\bf z'})\Vert^2\leq 3(K^2+1)( \vert x_i-x'_i\vert^2 +\vert y_i-y'_i\vert^2 )+3K^2 (\vert x_{i+1}-x'_{i+1}\vert^2+\vert y_{i-1}-y'_{i-1}\vert^2) $$ and  that
$$\Vert \xi({\bf z})-\xi({\bf z'})\Vert^2\leq (6K^2+3)\Vert {\bf z}-{\bf z'}\Vert^2.$$ \hfill$\Box$

\bigskip

One deduces that the associated differential system 
$$\cases{ \dot x_i= y_i-g'_{i-1}(x_{i},y_{i-1}),\cr \dot y_i= x_i-g_i(x_{i+1},y_i),\cr}$$
defines a flow on $E$. We will denote by ${\bf z}^t$ the image at time $t$ of a point 
${\bf z}\in E$ by this flow, and more generally by $Z^t$ the image of a subset $Z\subset E$.

\bigskip
\noindent{\sc Lemma 2:}\enskip\enskip {\it For every $({\bf z}, {\bf z}')\in E$ and every $t\in\R$, one has 
$$\e^{-A\vert t\vert}\Vert{\bf z}-{\bf z}'\Vert \leq\Vert {\bf z}^t-{\bf z}'{}^t\Vert \leq  e^{A\vert t\vert}\Vert{\bf z}-{\bf z}'\Vert $$and
$$e^{-A\vert t\vert}\Vert \xi({\bf z})\Vert \leq\Vert \xi({\bf z}^t)\Vert \leq  e^{A\vert t\vert}\Vert \xi({\bf z})\Vert 
.$$}

\bigskip
\noindent{\it Proof.} \enskip Let us begin with the first double inequality. For every $({\bf z}, {\bf z}')\in E$ and $t\geq 0$, one has
$$\eqalign {\Vert {\bf z}^t-{\bf z}'{}^t\Vert &\leq \Vert{\bf z}-{\bf z}'\Vert + \Vert ({\bf z}^t-{\bf z}'{}^t)-({\bf z}-{\bf z}')\Vert\cr
&= \Vert{\bf z}-{\bf z}'\Vert +\left\Vert\int_0^t \xi({\bf z}^s)- \xi({\bf z}'{}^s) \, ds\right\Vert\cr
&\leq \Vert{\bf z}-{\bf z}'\Vert + A\int_0^t \Vert{\bf z}^s- {\bf z}'{}^s\Vert \, ds \cr}$$
which implies by Gronwall's Lemma that $\Vert {\bf z}^t-{\bf z}'{}^t\Vert\leq  e^{At}\Vert{\bf z}-{\bf z}'\Vert $. The  inequality on the right, for $t\leq0$ can be proven similarly and the inequality on the left can be deduced immediately from the one on the right. 

Writing this double inequality with ${\bf z}'={\bf z}^s$, dividing by $s$ and letting $s$ tend to $0$ permits to obtain the second double inequality.
\hfill$\Box$

\bigskip
For every $i\in\Z$, define the maps
 $$\eqalign{Q_{i}: E&\to\R^2,\cr
 {\bf z}&\mapsto 
(g_{i}(x_{i+1},y_{i}), y_i),\cr}$$ and
$$\eqalign{Q'_i:E&\to \R^2,\cr
{\bf z}&\mapsto
(x_{i},g'_{i-1}(x_{i},y_{i-1})),\cr}
$$
and observe that $f_i\circ Q_i({\bf z})=Q'_{i+1}({\bf z}).$ Observe also that 
$$\xi_i=J\circ (Q'_i-Q_i),$$ where
$J(x,y)=(-y,x)$. In particular ${\bf z}$ is a singularity of $\xi$ if and only if $Q_i({\bf z})=Q'_{i}({\bf z})$ for every $i\in\Z$. One deduces that $Q_{1}$ induces a bijection between the set of
singularities of
$\xi
$ and the set of fixed points of $f^q$. Indeed, if ${\bf z}$ is a
singularity of $\xi$, we have
$$f^q\circ Q_1({\bf z})=f_{mq}\circ\dots\circ f_1\circ
Q_1({\bf z})=Q_{mq+1}({\bf z})=Q_1({\bf z}).$$
The set of singularities consists of the constant sequence ${\bf 0}=(0)_{i\in\Z}$, whose image by $Q_1$ is the common fixed point $0$ of all $f_i$, and of finitely many smooth closed curves $(\Sigma_p)_{p\in (q\alpha, q\beta)\cap {\Z}}$, each curve $\Sigma_p$ being sent homeomorphically onto $S_{p/q}$ by $Q_1$ (and in fact by each $Q_i$ or $Q'_i$). Observe that $\xi$ is $C^1$ in a neighborhood of ${\bf 0}$ and $C^{\infty}$ in a neighborhood of $\Sigma_p$ because each $f_i$ is a $C^1$ diffeomorphism in a neighborhood of $0$ and a $C^{\infty}$ diffeomorphism in a neighborhood of $S_{p/q}$.

\medskip

Observe that $\xi$ is the gradient
flow of the function
$${\bf h}:{\bf z}\mapsto\sum_{0<i\leq mq} x_iy_i-h_{i-1}(x_i,y_{i-1})$$ and that ${\bf h}$ is invariant by $\varphi$. The {\it energy} of an orbit $({\bf z}^t)_{t\in\R}$ is
$$ \int_{-\infty}^{+\infty} \Vert \xi({\bf z}^t)\Vert^2 \, dt= \lim_{t\to+\infty}  {\bf h}({\bf z^t})- \lim_{t\to-\infty}  {\bf h}({\bf z^t}).$$
Using the inequality $$e^{-A\vert t\vert}\Vert \xi({\bf z})\Vert\leq \Vert \xi({\bf z}^t)\Vert $$ in the integral, one gets$$\Vert \xi({\bf z})\Vert^2\leq  A\int_{-\infty}^{+\infty} \Vert \xi({\bf z}^t)\Vert^2 \, dt = A\left(\lim_{t\to+\infty}  {\bf h}({\bf z^t})- \lim_{t\to-\infty}  {\bf h}({\bf z^t})\right).$$

\bigskip
\noindent{\sc Lemma 3:}\enskip\enskip {\it For every ${\bf z}\in\Gamma_p$, one has
$${\bf h}({\bf z})-{\bf h}({\bf 0})=\pi(p-q \alpha) \left( 1 + (p/q -\alpha) +{(p/q-\alpha)^2\over 3}\right).$$}

\bigskip
\noindent{\it Proof.} \enskip The quantity ${\bf h}({\bf z})-{\bf h}({\bf 0})$ is the difference of {\it action} between the two corresponding fixed points of $f^q$. It measures the area diplaced by an arc joining $0$ to $Q_1(z)$ along an isotopy of $\R^2$ that fixes $0$ and $Q_1(z)$. In our case, it is independent of ${\bf z}$, the set $\Sigma_p$ being contained in the critical set of ${\bf h}$, and its dynamical meaning is the following: it is equal to the opposite of the rotation number of the Lebesgue measure in the annulus  $D_{p/q}\setminus\{0\}$ defined by the map $f^q$ for the lift to the universal covering space that fixes all the points of the boundary line. It is easy to compute. The rotation number of the Lebesgue measure restricted to ${\bf \D}\setminus\{0\}$ is equal to $\pi(q\alpha-p)$, therefore: 
 $$\eqalign{{\bf h}({\bf z})-{\bf h}({\bf 0})&=\pi(p-q \alpha)-\int_1^{1+p/q-\alpha} 2\pi r((r+\alpha-1)q-p)\,dr\cr
 &= \pi(p-q \alpha) \left( 1 + (p/q -\alpha) +{(p/q-\alpha)^2\over 3}\right).\cr}$$
\hfill$\Box$

\bigskip
\bigskip

\noindent{\large {\bf 3. The main proposition and its consequences}}

\bigskip
In this section we will give a proof of Theorem 1'. It will  follow from Proposition 2, whose proof is postponed in Section 5. The arguments are very close to the ones given by Bramham in \cite{Br1}: we replace pseudoholomorphic curves by orbits of the gradient flow but the spirit of the proof is the same.

\bigskip
\noindent{\sc Proposition 2:}\enskip\enskip {\it The curve $\Sigma_p$ bounds a topological disk $\Delta_{p}$ that satisfies the following:

\smallskip 
\noindent{\bf i)}\enskip $\Delta_{p}$ contains $0$;

\smallskip 
\noindent{\bf ii)}\enskip $\Delta_{p}$ is invariant by $\varphi$;

\smallskip 
\noindent{\bf iii)}\enskip each projection ${\bf z}\mapsto (x_i,y_{i-1})$, $i\in\Z$, is one to one on $\Delta_{p}$;

\smallskip 
\noindent{\bf iv)}\enskip each projection ${\bf z}\mapsto (x_i,y_i)$, $i\in\Z$,  is one to one on $\Delta_{p}$;

\smallskip 
\noindent{\bf v)}\enskip $\Delta_{p}$ is invariant by the flow;

\smallskip 
\noindent{\bf vi)}\enskip for every ${\bf z}\in \Delta_{p}\setminus (\{{\bf 0}\}\cup\Sigma_p)$, one has $\lim_{t\to-\infty} {\bf z}^t={\bf 0}$ and $\lim_{t\to+\infty} d({\bf z}^t, \Sigma_p)=0$. 

}

\bigskip
Let us explain now why this proposition implies Theorem 1'.

\bigskip
\noindent{\it Proof of Theorem 1'.} \enskip  
The assertion {\bf iii)} tells us that the maps $Q_i\vert_{\Delta_{p}}$ and $Q'_i\vert_{\Delta_{p}}$, $i\in\Z$, induce homeomorphisms from $\Delta_{p}$ to $D_{p/q}$.
One gets a family of homeomorphisms $(\widehat f_i)_{i\in\Z}$ of $D_{p/q}$ by writing:
$$\widehat f_i=(Q_{i+1}\vert_{\Delta_{p}})\circ (Q_i\vert _{\Delta_{p}})^{-1}.$$This family is $m$ periodic because ${\Delta_{p}}$ is invariant by $\varphi$. Morevover  $\widehat f =\widehat f_m\circ\dots\circ \widehat f_1$ is $q$ periodic because
$$\widehat f =(Q_{m+1}\vert_{\Delta_{p}})\circ (Q_{1}\vert_{\Delta_{p}})^{-1}=(Q_{1}\vert_{\Delta_{p}})\circ (\varphi\vert_{\Delta_q})\circ (Q_{1}\vert_{\Delta_{p}})^{-1}.$$
Observe now that $$f_i\vert_{D_{p/q}}=(Q'_{i+1}\vert_{\Delta_{p}})\circ (Q_i\vert {\Delta_{p}})^{-1}$$
and that 
$$\widehat f_i-f_i\vert_{D_{p/q}}= J \circ \xi_{i+1}\circ (Q_i\vert {\Delta_{p}})^{-1}.$$
Observe also that $\widehat f_i$ fixes $0$ and coincides with $f_i$ on $S_{p/q}$.

\medskip
On the disk $D_{p/q}$, one can write
$$\eqalign{\widehat f - f=\enskip\enskip &\widehat f_m\circ\widehat f_{m-1}\circ\dots\circ \widehat f_1\enskip - \enskip f_m\circ\widehat f_{m-1}\circ\dots\circ \widehat f_1 \cr
 +\enskip & f_m\circ\widehat f_{m-1}\circ\dots\circ \widehat f_1 \enskip - \enskip f_m\circ f_{m-1}\circ\dots\circ \widehat f_1\cr
 +\enskip &\dots\dots\dots\cr
 + \enskip&f_m\circ\dots\circ f_{2}\circ \widehat f_1\enskip - \enskip f_m\circ\dots\circ f_{2}\circ  f_1.\cr}$$
By definition of $K$ Lipschitz untwisted homeomorphisms, one deduces that
$$\sup_{z\in D_{p/q}} \vert \widehat f(z)-f(z)\vert\leq (1+K+ \dots+ K^{m-1})A^{1/2} C(p,q)^{1/2},$$
where
$$C(p,q)= \pi(p-q \alpha) \left( 1 + (p/q -\alpha) +{(p/q-\alpha)^2\over 3}\right).$$

Let us consider now the homothety $H$ of ratio $1+p/q-\alpha$ and set $\check f=H^{-1}\circ \widehat f\circ H$. On the disk $\D$ one can write
$$\check f- f= (H^{-1}\circ \widehat f\circ H-H^{-1}\circ f\circ H)+(H^{-1}\circ f\circ H- f\circ H)+(f\circ H-f)$$
 and gets
 $$\sup_{z\in \D} \vert \check f(z)-f(z)\vert\leq (1+K+ \dots+ K^{m-1})A^{1/2} C(p,q)^{1/2}+ (1+K^m) (p/q-\alpha).$$
We will get the same upper bound for  $\sup_{z\in \D} \vert \check f^{-1}(z)-f^{-1}(z)\vert$. One can choose $p$ and $q$, with $p-q\alpha$ arbitrarily small, which means that the quantity on the right side itself can be chosen arbitrarily small. \hfill$\Box$

\bigskip
\noindent{\bf Remark.} \enskip  Keeping the notations above, one gets
$$\eqalign{\mathrm{Id}- f^q=\enskip\enskip &\widehat f_{mq}\circ\widehat f_{mq-1}\circ\dots\circ \widehat f_1\enskip - \enskip f_{m}\circ\widehat f_{mq-1}\circ\dots\circ \widehat f_1 \cr
 +\enskip & f_{mq}\circ\widehat f_{mq-1}\circ\dots\circ \widehat f_1 \enskip - \enskip f_{mq}\circ f_{mq-1}\circ\dots\circ \widehat f_1\cr
 +\enskip &\dots\dots\dots\cr
 + \enskip&f_{mq}\circ\dots\circ f_{2}\circ \widehat f_1\enskip - \enskip f_{mq}\circ\dots\circ f_{2}\circ  f_1.\cr}$$
which implies that
$$\sup_{z\in D_{p/q}} \vert z-f^q(z)\vert\leq (1+K+ \dots+ K^{mq-1})A^{1/2} C(p,q)^{1/2},$$
and that  $$\sup_{z\in \D} \vert z-f^q(z)\vert\leq (1+K+ \dots+ K^{mq-1})A^{1/2} C(p,q)^{1/2}+ (1+K^{mq}) (p/q-\widetilde \alpha).$$
One gets a similar inequality for  $\sup_{z\in \D} \vert z-f^{-q}(z)\vert$.  Extending in a similar way our original diffeomorphism with the help of a negative polar twist map will give us a similar inequality for couples $(p,q)$ such that $p/q<\alpha$.  If $\alpha$ satisfies the following super Liouville condition: for every $\mu\in(0,1)$, there exists two sequences of integers $(q_n)_{n\geq 0}$ and $(p_n)_{n\geq 0}$, with $q_n>0$, such that $\vert q_n\alpha-p_n\vert \leq \mu^{q_n}$, then there exists a sequence  $(r_n)_{n\geq 0}$ such that $(f^{r_n})_{n\geq 0}$ converges to the identity on $\D$ for the $C^0$ topology. One says that $f$ is {\it $C^0$ rigid}. This is a $C^1$ version, but with the additional assumption of being $C^1$ conjugate to the rotation, of the following recent result of Bramham \cite{Br2}:

\bigskip
\noindent{\sc Theorem 2:}\enskip\enskip {\it Every $C^{\infty}$ irrational pseudo rotation $f$ of rotation number $\overline\alpha$ is $C^0$ rigid if it satisfies the following super Liouville condition: for every $\mu\in(0,1)$, there exists a sequence of integers $(q_n)_{n\geq 0}$ such that $d(q_n\overline \alpha, 0) \leq \mu^{q_n}$}

\bigskip
\bigskip

\noindent{\large {\bf 4. Canonical dominated structure for the gradient flow}}

\bigskip
In this section, we will do a deeper study of the vector field $\xi$. The fact that $\xi$ is a gradient flow has no importance here. What is crucial is the fact that $\xi$ is {\it tridiagonal and monotonically symmetric} (we will explain the meaning). We refer to \cite {Lec1} or \cite{Lec2}  for detailed proofs. In what follows, the function $\mathrm{sign}$ assigns $+1$ to a positive number and $-1$ to a negative number.

Let us consider the set
$$V=\{{\bf z}\in E\enskip\vert
x_i\not=0\enskip \mathrm{and}\enskip  y_i\not=0\enskip{\rm for\enskip all}\enskip i\in\Z\}$$ and the function $L$ on $V$ defined by the formula
$$\eqalign{L({\bf z})&={1\over 4}\sum_{0<i\leq mq} \mathrm{sign}( x_i) \left( \mathrm{sign}(y_i)-\mathrm{sign} (y_{i-1})\right),\cr&= {1\over 4}\sum_{0<i\leq mq}  \mathrm{sign}( y_i) \left( \mathrm{sign}(x_i)-\mathrm{sign} (x_{i+1})\right).\cr}$$It extends
continuously on the open set 
$$W =\{{\bf z}\in E\enskip\vert\enskip x_i=0\Rightarrow y_{i-1}y_{i}>0, \enskip  y_i=0\Rightarrow x_{i}x_{i+1}>0\}.$$
Let us explain first what is the meaning of $L$. For every $z\in E$, one can define a loop $\gamma_{\bf z}: [0,2mq]\to\R^2$ by writing:
$$\gamma_{\bf z}(t)=\cases {\left((1+2i-t) x_i+ (t-2i) x_{i+1}, y_i\right) & if $t\in[2i,2i+1]$,\cr
\left(x_{i+1}, (2+2i-t) y_i+ (t-2i-1) y_{i+1}\right) & if $t\in[2i+1,2i+2]$.\cr}$$
The fact that ${\bf z}$ belongs to $W$ means that the image of this loop does not meet $0$.  Write $0x^+$, $0x^-$, $0y^+$, $0y^-$ for the half lines generated by the vectors $(1,0)$, $(-1,0)$, $(0,1)$, $(0,-1)$ respectively. The formulas given above tells us that
$$L({\bf z})=  {1\over 2} (0x^+\wedge\gamma_z +0x^-\wedge\gamma_z)={1\over 2} (0y^+\wedge\gamma_z +0y^-\wedge\gamma_z) , $$  where $\wedge$ means the algebraic intersection number. The integer $L(z)$ is nothing but the indice of the loop  $\gamma_{\bf z}$ relatively to $0$. In particular, $L$ is integer valued and takes its values in $\{-[nq/2],\dots, [nq/2]\}$.

\medskip 
Let us state the fundamental result, whose proof is postponed in the appendix:

\bigskip
\noindent{\sc Proposition 3~:}\enskip\enskip{\it If
${\bf z}$, 
${\bf z}$ are two distinct points of $E$ such that ${\bf z}'-{\bf z}\not\in W$, then there exists $\varepsilon >0$ such
that
${\bf z'}^{t}-{\bf z}^{t}\in W$ if  $0< \left \vert t\right 
\vert \leq \varepsilon $ and one has
$$L({\bf z'}^{\varepsilon} -{\bf z}^{\varepsilon 
})>L({\bf z'}^{-\varepsilon } -{\bf z}^{-\varepsilon })\enskip.$$ }

\bigskip
For every $p\in \{-[nq/2],\dots, [nq/2]\}$, let us write $W_p=\{{\bf z}\in W\,\vert\, L({\bf z})=p\}$ and  define
 $$\enskip W_p^{+}= \mathrm{Int}\left(\mathrm{Cl}\left(\bigcup_{p'\geq p} W_{p'}\right)\right) , \enskip W_p^{-}= \mathrm{Int}\left(\mathrm{Cl}\left(\bigcup_{p'\leq p} W_{p'}\right)\right) .$$
 
Similarly write 
$$\eqalign{{\cal W}_p&=\{({\bf z}, {\bf z}')\in E\times E\enskip\vert\enskip {\bf z}'-{\bf z}\in W_p\},\cr
{\cal W}_p^{+}&=\{({\bf z}, {\bf z}')\in E\times E\enskip\vert\enskip {\bf z}'-{\bf z}\in W_p^{+}\},\cr
{\cal W}_p^{-}&=\{({\bf z}, {\bf z}')\in E\times E\enskip\vert\enskip {\bf z}'-{\bf z}\in W_p^{-}\}.\cr}$$

This result asserts that there exists a canonical filtration on the product flow defined on $E\times E\setminus \mathrm{diag}$, where $\mathrm{diag}=\{({\bf z}, {\bf z}')\in E\times E\,\vert\, {\bf z}={\bf z}'\}$.  Each set  ${\cal W}_p^{+}$ is an attracting set of the product flow on $E\times E\setminus \mathrm{diag}$ and each set ${\cal W}_p^{-}$ a repulsing set. More precisely, if $({\bf z}, {\bf z}')\in\mathrm{Cl}({\cal W}_p^{+})\setminus\mathrm{diag}$, then  $({\bf z}^t, {\bf z}'{}^t)\in{\cal W}_p^{+}$ for every $t>0$; if $({\bf z}, {\bf z}')\in\mathrm{Cl}({\cal W}_p^{-})\setminus\mathrm{diag}$, then  $({\bf z}^t, {\bf z}'{}^t)\in{\cal W}_p^{-}$ for every $t<0$. Consequently, the boundary of ${\cal W}_p^{+}$ and ${\cal W}_p^{-}$ in $E\times E\setminus\mathrm{diag}$ are 1 codimensional topological submanifolds. 
 
 \medskip
 In particular, if ${\bf z}$ and ${\bf z'}$ are two singularities, then ${\bf z}'-{\bf z}\in W$ and $L({\bf z}'-{\bf z})$ is well defined. Let us explain the meaning of this integer. Recall that to each map $f_:(x,y)\to(X,Y)$ is naturally  associated an isotopy $(f_{i,t})_{t\in[0,2]}$ starting from the identity defined as follows $$f_{i,t}(x,y)= \cases {((1-t) x+t X, y) & if $t\in[0,1]$,\cr
 (X, (2-t)y+ (t-1) Y) & if $t\in[1,2]$,\cr}$$ and an isotopy $(f^{[q]}_t)_{t\in[0,2mq]}$ joining the identity to $f^q$, where
 $$f^{[q]}_t= f_{1+i,t-2i}\circ f_{i}\circ \dots \circ f_1 \enskip\mathrm{if}\enskip t\in [2i, 2i+2].$$
The integer $L({\bf z}'-{\bf z})$ is equal to the {\it linking number} of the two corresponding fixed points of $f^q$ for this natural isotopy naturally defined by the decomposition of $f$. In particular, if  $p/q\in(\alpha , \beta)$, then $({\bf 0}, {\bf z})\in {\cal W}_p$ for every $z\in \Sigma_p$ and $({\bf z}, {\bf z}')\in {\cal W}_{p}$ for every ${\bf z}$ and ${\bf z}'$ in $\Sigma_p$.

 \medskip
 The reason why Proposition 3 is true is the fact that the vector field is tridiagonal and monotonically symmetric. Writing the coordinates in the following order
 $$ \dots, y_{i-1}, x_i, y_i, x_{i+1}, \dots $$
 the corresponding coordinate of $\xi$ depends only on this coordinate and its two neighbours. Moreover it depends monotonically of each of the neighbouring  coordinates and for two neighbouring coordinates, the ``cross monotonicities" are the same. In our example, $\dot x_i$ depends only on $x_i$, $y_{i-1}$ and $y_{i}$, is a decreasing function of $y_{i-1}$ and an increasing function of $y_{i}$ whereas $\dot y_i$ depends only on $y_i$, $x_i$ and $x_{i+1}$, is an increasing function of $x_{i}$ and a decreasing function of $x_{i+1}$.  To every tridiagonal and monotonically symmetric vector field is associated a natural function $L$ satisfying Proposition 3. An important case is the linear case. Suppose that $\xi_*$ is a linear tridiagonal and monotonically symmetric vector field on our space $E$. We obtain a dominated splitting (that has been known for a long time, see \cite{M} for example): there exists  a linear decomposition
 $$E=\bigoplus_{p\in\{-[nq/2],\dots, [nq/2]\}} E_{p}$$ in invariant subspaces where $$E_{p}\setminus\{0\}=\left\{{\bf z}\in E\, \vert \,e^{t\xi_*}({\bf z})\in W_p\enskip\mathrm{for\enskip all} \enskip t\in\R\right\},$$and the real parts of the eigenvalues of $\xi_*\vert_{E_{*,p'}}$ are larger than the real parts of the eigenvalues of $\xi_*\vert_{E_{p}}$, if $p'<p$. 
Moreover,  the spaces
$$E_{p}^{+}= \bigoplus_{p'\geq p} E_{p'},\enskip E_{p}^{-}= \bigoplus_{p'\leq p} E_{p'}$$satisfy
$$E_{p}^{+}\setminus\{0\}=\left\{{\bf z}\in E\, \vert \,e^{t\xi_*}({\bf z})\in W_p^{+}\enskip\mathrm{for\enskip all} \enskip t\in\R\right\}$$and
$$E_{p}^{-}\setminus\{0\}=\left\{{\bf z}\in E\, \vert \,e^{t\xi_*}({\bf z})\in W_p^{-}\enskip\mathrm{for\enskip all} \enskip t\in\R\right\}.$$  In the case of a tridiagonal and monotonically symmetric $C^1$ vector field, with non zero cross derivatives, one gets such a decomposition of the tangent space at every singularity, for the linearized flow. The proof of Proposition 3, in the $C^1$ case is given in \cite{Lec1}.  Starting with two distinct points ${\bf z}$ and ${\bf z}'$ such that such that ${\bf z}'-{\bf z}\not\in W$,
polynomial approximations obtained by successive integrations permit to determine the sign of the coordinates of  $z^t-z'{}^t$, for small values of $t$. A more precise study, replacing polynomial approximations by explicit inequalities is given in \cite{Lec2} (Lemma 2.5.1) and permits to extend Proposition 3 to the compactification of $E\times E\setminus \mathrm{diag}$ obtained by blowing up the diagonal and then to get a similar result for the linearized vector field on the tangent bundle. One proves in that way the existence of a global dynamically coherent dominated splitting (the definition is recalled in the next section). The proof of Proposition 3 that uses Lemma 2.5.1 of \cite{Lec2} extends word to word to our Lipschitz case. However the following proposition, that will be needed in the next section, necessitates an extension of this lemma  to the case where vanishing conditions on the coordinates of $z-z'$, are replaced by smallness conditions. This will be the object of Lemma 6 in  the appendix.

\bigskip
\bigskip
Let us define
$$E^i=\{{\bf z}\in E\,\vert\enskip i'\not=i\Rightarrow x_{i'}=0\enskip\mathrm{and}\enskip i'\not=i-1\Rightarrow y_{i'}=0\}$$
and 
$$E'{}^i=\{{\bf z}\in E\,\vert\enskip i'\not=i\Rightarrow x_{i'}=0\enskip\mathrm{and}\enskip i'\not=i\Rightarrow y_{i'}=0\}.$$
Write $\pi^{i}:E\to E^i$ and $\pi'{}^{i}:E\to E'{}^i$ for the orthogonal projections on $E^i$ and $E'{}^i$ respectively, write $\pi^{i}{}^{\perp}:E\to E^i{}^{\perp}$ and $\pi'{}^{i}{}^{\perp}:E\to E'{}^i{}^{\perp}$ for the orthogonal projections on $E^i{}^{\perp}$ and $E'{}^i{}^{\perp}$.

\bigskip
\noindent{\sc Proposition 4:}\enskip\enskip{\it For every $t>0$, there exists a constant $N_t>0$ such that
for every points
${\bf z}$, 
${\bf z}'$ in $E$ satisfying ${\bf z}'^{s}-{\bf z}^s\in W$ if $s\in[-t,t]$ and for every $i\in\Z$, one has
$$\Vert\pi^i{}^{\perp}({\bf z'}-{\bf z})\Vert\leq N_t \,\Vert\pi^i({\bf z'}-{\bf z})\Vert$$ 
and
$$\Vert\pi'{}^i{}^{\perp}({\bf z'}-{\bf z})\Vert\leq N_t \,\Vert\pi'{}^i({\bf z'}-{\bf z})\Vert.$$ }

\bigskip
\bigskip

\noindent{\large {\bf 5. Proof of Proposition 2}}

\bigskip
The goal of this section is to prove Proposition 2. We will begin with a preliminary result (Proposition 5) that states the existence of a topological plane $\Pi_p$, containing ${\bf 0}$ and $\Sigma_p$, invariant by $\varphi$, such that for every $t\in\R$ and every $i\in\Z$, the projections ${\bf z}\mapsto (x_i,y_{i-1})$, and ${\bf z}\mapsto (x_i,y_{i})$ send homeomorphically $(\Pi_p)^t$ onto $\R^2$. This last fact will be a consequence of the inclusion $$\left((\Pi_p)^t\times (\Pi_p)^t\right)\setminus\mathrm{diag}\subset{\cal W}_p.$$  The disk $\Delta_{p}\subset \Pi_p$ bounded by $\Sigma_p$ satisfies the assertions {\bf i)} to {\bf iv)} of Proposition 2. In the second part of the section, we will prove that it is invariant by the flow (assertion {\bf v)}) and that for every ${\bf z}\in \Delta_{p}\setminus (\{{\bf 0}\}\cup\Sigma_p)$, one has $\lim_{t\to-\infty} {\bf z}^t={\bf 0}$ and $\lim_{t\to+\infty} d({\bf z}^t, \Sigma_p)=0$ (assertion {\bf vi)}).

As told in the previous section, a tridiagonal and monotonically symmetric $C^1$ vector field on $E$, with non zero cross derivatives, admits a dominated splitting, which means a decomposition of the tangent bundle
 $$TE=\bigoplus_{p\in\{-[mq/2],\dots, [mq/2]\}} E_{p}({\bf z}),$$
invariant by the linearized flow, with relative expanding properties, and this splitting is dynamically coherent, which means that every field
$$\bigoplus_{p_0\leq p\leq p_1} E_{p}({\bf z}),$$is integrable. This is a consequence of Proposition 3 stated in \cite{Lec2}.  The coherency is obtained via graph transformations. One can integrate the fields 
$$\bigoplus_{ p_0\leq p} E_{p}({\bf z}),\enskip \bigoplus_{p\leq p_1} E_{p}({\bf z}),$$ and then take the intersection of the integral manifolds. In particular a plane $\Pi_p$ tangent to the bundle $E_p({\bf z})$ is characterized by the property
$$\left((\Pi_p)^t\times (\Pi_p)^t\right)\setminus\mathrm{diag}\subset{\cal W}_p\enskip\mathrm{for \enskip all} \enskip t\in\R.$$ 
Here the vector field is no more $C^1$, the decomposition 
$$TE=\bigoplus_{p\in\{-[mq/2],\dots, [mq/2]\}} E_{p}({\bf z})$$
does not exist but fortunately the graph transformation exists: there exist topological manifolds $\Gamma^+_p$ and $\Gamma^-_p$ satisfying 
$$\mathrm{dim }(\Gamma^+_p)+ \mathrm{dim }(\Gamma^-_p)=mq+2$$
and
$$\left((\Gamma^+_p)^t\times (\Gamma^+_p)^t\right)\setminus\mathrm{diag}\subset{\cal W}^{+}_p, \enskip\left((\Gamma^-_p)^t\times (\Gamma^-_p)^t\right)\setminus\mathrm{diag}\subset{\cal W}^{-}_p$$
for all $t\in\R$. 
In case where the vector field is $C^1$, the manifolds $\Gamma^+_p$ and $\Gamma^-_p$ are  $C^1$ and the fact that $\Pi_p=\Gamma^+_p\cap\Gamma^-_p$
 is a $C^1$ plane follows almost immediately from the Implicit Function Theorem. In our situation, to prove that $\Pi_p=\Gamma_p^+\cap\Gamma_p^-$ is a plane, we will need a Lipschitz Implicit Function Theorem which means that some Lipschitz conditions about the manifolds $\Gamma^+_p$, $\Gamma^-_p$ and $\Pi_p$ must be satisfied. This is the reason one needs Proposition 4.  

\bigskip

\noindent{\sc Proposition 5:}\enskip\enskip {\it  For every $p\in(q\alpha, q\beta)\cap \Z$, there exists a set $\Pi_p\subset E$, image of a proper topological embedding of $\R^2$, such that:

\smallskip 
\noindent{\bf i)}\enskip $\Pi_p$ contains $\{{\bf 0}\}\cup\Sigma_p$;

\smallskip 
\noindent{\bf ii)}\enskip $\Pi_p$ is invariant by $\varphi$;

\smallskip 
\noindent{\bf iii)}\enskip $((\Pi_p)^t\times (\Pi_p)^t)\setminus\mathrm{diag}\subset{\cal W}_p$, for every $t\in\R$.
}

\bigskip
\noindent{\it Proof.} \enskip As explained in the Remark 2 at the end of Section 1, the rotation $f_{*,i}$ that coincides with $f_i$ on the circle $S_{p/q}$ is an area preserving $K$ Lipschitz untwisted homeomorphism. One gets a decomposition $f_*=f_{*,m}\circ\dots\circ f_{*,1}$ of the rotation of angle $2\pi p/q$. To this decomposition is associated a linear vector field $\xi_*$ on $E$ which is the gradient of a quadratic form. Its kernel being homeomorphic to the fixed point set of $(f_*)^q$, is a plane. Since $\xi_*$ coincides with $\xi$ on $\{0\}\cup\Sigma_p$, its kernel contains this set: it is the plane generated by $\{0\}\cup\Sigma_p$. Applying what has been said in Section 4, one knows that there exists a linear (and orthogonal) decomposition
 $$E=\bigoplus_{p'\in\{-[nq/2],\dots, [nq/2]\}} E_{p'},$$
where
$$E_{p'}\setminus\{0\}=\left\{{\bf z}\in E\, \vert \,e^{t\xi_*}({\bf z})\in W_{p'}\enskip\mathrm{for\enskip all} \enskip t\in\R\right\},$$ and that the spaces
$$E_{p'}^{+}= \bigoplus_{p''\geq p'} E_{p''},\enskip E_{p'}^{-}= \bigoplus_{p''\leq p'} E_{p''},$$satisfy
$$E_{p'}^{+}\setminus\{0\}=\left\{{\bf z}\in E\, \vert \,e^{t\xi_*}({\bf z})\in W_{p'}^{+}\enskip\mathrm{for\enskip all} \enskip t\in\R\right\}$$and
$$E_{p'}^{-}\setminus\{0\}=\left\{{\bf z}\in E\, \vert \,e^{t\xi_*}({\bf z})\in W_{p'}^{-}\enskip\mathrm{for\enskip all} \enskip t\in\R\right\}.$$ Observe that $E_{p}$ is the kernel of $\xi_*$ and that
$E_{p}$, $E_{p}^{+}$ and $E_{p}^{-}$ are invariant by $\varphi$ because $\xi_*$ is invariant by $\varphi$.
 Write $$\pi_{p}^{+}: E\to E_{p}^{+}, \enskip \pi_{p}^{-}: E\to E_{p}^{-}$$ for the orthogonal projections. Every vector ``sufficiently close'' to  $E_{p}^{+}$ or $E_{p}^{-}$ must belong to $W_{p}^{+}$ or $W_{p}^{-}$ respectively. This implies that there exists a constant $M\geq 1$ such that:
 $$({\bf z},{\bf z}')\in {\cal W}_p^{+}\Rightarrow \Vert \pi_{p-1}^{-}({\bf z})-\pi_{p-1}^{-}({\bf z}')\Vert \leq M\Vert \pi_{p}^{+}({\bf z})-\pi_{p}^{+}({\bf z}')\Vert,$$
 and
 $$({\bf z},{\bf z}')\in {\cal W}_p^{-}\Rightarrow \Vert \pi_{p+1}^{+}({\bf z})-\pi_{p+1}^{+}({\bf z}')\Vert \leq M\Vert \pi_{p}^{-}({\bf z})-\pi_{p}^{-}({\bf z}')\Vert.$$

\medskip

Identifying the space $E=E_{p}^{+} \oplus E_{p-1}^{-}$ with the product $E_{p}^{+}\times E_{p-1}^{-}$, we write $\Gamma_{\psi}\subset E$ for the graph of a fonction $\psi: E_{p}^{+} \to  E_{p-1}^{-}$. We define
 $${\cal G}_{p}^{+}= \left\{\psi: E_{p}^{+} \to  E_{p-1}^{-}\enskip\vert\enskip(\Gamma_{\psi}\times \Gamma_{\psi})\setminus\mathrm{diag}\subset{\cal W}_p^{+}\right\}$$ and
 $$\overline{\cal G}_{p}^{+}= \left\{\psi: E_{p}^{+} \to  E_{p-1}^{-}\enskip\vert\enskip(\Gamma_{\psi}\times \Gamma_{\psi})\setminus\mathrm{diag}\subset\mathrm{Cl}({\cal W}_p^{+})\right\}.$$ 
 
Note that $\overline{\cal G}_{p}^{+}$ is closed for the compact-open topology, note also that every function in $\overline{\cal G}_{p}^{+}$ is  $M$ Lipschitz.

\bigskip
\noindent{\sc Lemma 4:}\enskip\enskip {\it The vector field $\xi$ induces a positive semi-flow 
$(t,\psi)\mapsto \psi^t$ on $\overline{\cal G}_{p}^{+}$, such that $\Gamma_{\psi^t}=(\Gamma_{\psi})^t$. Moreover, one has $\psi^t\in{\cal G}_{p}^{+}$ for every $\psi\in \overline{\cal G}_{p}^{+}$ and every $t>0$.}

\bigskip
\noindent{\it Proof.} \enskip By Proposition 3, for every $\psi\in \overline{\cal G}_{p}^{+}$ and  every $t\geq 0$, one has 
$$\left((\Gamma_{\psi})^t\times (\Gamma_{\psi})^t\right)\setminus\mathrm{diag}\subset\mathrm{Cl}({\cal W}_p^{+}).$$ Consequently $(\Gamma_{\psi})^t$ projects injectively into $E_p^{+}$. In fact it projects surjectively. Indeed, the map
$${\bf z}\mapsto \pi_{p} (({\bf z}+ \psi({\bf z}))^t)$$ is more than an injective and continuous transformation of  $E_{p}^{+}$. By Lemma 2, it is  $Me^{At}$ bi-Lipshitz. In particular it is a homeomorphism of $E_{p}^{+}$. This means that $(\Gamma_{\psi})^t$ is the graph  a continuous fonction $\psi^t \in{\cal G}_p^{+}$. The continuity of the map
$$(t,\psi)\mapsto \psi^t,$$ when $\overline{\cal G}_{p}^{+}$ is munished with the compact-open topology, follows easily. The fact that $\psi^t$ belongs to ${\cal G}_p^{+}$ for every $t>0$ is an immediate consequence of Proposition 3.\hfill$\Box$

\bigskip

Similarly, one can identify $E$ with the product $E_{p}^{-}\times E_{p+1}^{+}$, write $\Gamma'_{\psi'}\subset E$ for the graph of a fonction $\psi': E_{p}^{-} \to  E_{p+1}^{+}$  and define $${\cal G}_p^{-}= \left\{\psi': E_{p}^{-} \to  E_{p+1}^{+}\enskip\vert \enskip (\Gamma'_{\psi'}\times\Gamma'_{\psi'})\setminus\mathrm{diag}\subset{\cal W}_p^{-}\right\}$$ and
 $$\overline{\cal G}_{p}^{-}= \left\{\psi: E_{p}^{-} \to  E_{p+1}^{+}\enskip\vert \enskip (\Gamma'_{\psi'}\times \Gamma'_{\psi'})\setminus\mathrm{diag}\subset\mathrm{Cl}({\cal W}_p^{-})\right\}.$$ 
We can define  a negative semi-flow 
$(t,\psi')\mapsto \psi'{}^t$ on $\overline{\cal G}_{p}^{-}$ such that $\Gamma'_{\psi'{}^t}=(\Gamma'_{\psi'})^t$ and we have $\psi'{}^t\in{\cal G}_{p}^{-}$ for every $\psi'\in \overline{\cal G}_{p}^{-}$ and every $t<0$.

\medskip
If ${\bf z}\in E_{p-1}^{-}$ denote by $\psi_{\bf z}\in {\cal G}_p^{+}$ the constant map equal to ${\bf z}$, and similarly if ${\bf z}'\in E_{p+1}^{+}$ denote by $\psi'_{{\bf z}'}\in {\cal G}_p^{-}$ the constant map equal to ${\bf z}'$. The graphs of these families of maps give us two transverse foliations. Let us study the time evolution of these foliations.
For every $t>0$, every ${\bf z}\in E_{p-1}^{-}$ and every  ${\bf z}'\in E_{p+1}^{+}$, we define the set 
$$\Pi_{{\bf z}, {\bf z}',t}=\Gamma_{\psi_{\bf z}^t}\cap \Gamma'_{\psi'{}_{{\bf z}'}^{-t}}.$$

\bigskip
\noindent{\sc Lemma 5:}\enskip\enskip {\it  Let us fix $i\in\Z$. The set $\Pi_{{\bf z}, {\bf z}',t} $ is the graph of a map $\theta_{i,t}: E^{i}\to E^{i}{}^{\perp}$ and the graph of a map $\theta'_{i,t}: E'{}^{i}\to E'{}^{i}{}^{\perp}$. Moreover, $\theta_{i,t}$ and $\theta'_{i,t}$ are $N_t$ Lipschitz, where $N_t$ is defined by Proposition 4.}

\bigskip
\noindent{\it Proof.} \enskip For every $s\in[-t,t]$, one has
 $$(\Pi_{{\bf z}, {\bf z}',t})^s\subset \Gamma_{\psi_{\bf z}^{t+s}}\cap \Gamma'_{\psi'{}_{{\bf z}'}^{-t+s}}$$which implies that 
$$\left((\Pi_{{\bf z}, {\bf z}',t})^s\times (\Pi_{{\bf z}, {\bf z}',t})^s\right)\setminus\mathrm{diag}\subset {\cal W}_p^{+}\cap {\cal W}_p^{-}={\cal W}_p.$$
Therefore $\Pi_{{\bf z}, {\bf z}',t}$ projects injectively on $E^i$ and $E'{}^i$. We want to prove that it projects surjectively. Fix ${\bf z_*}\in E^i$ and look at the map
$$\eqalign{ \Theta: E^i{}^{\perp}&\to E_{p-1}^{-}\times E_{p+1}^{+},\cr
{\bf z} &\mapsto \left (\pi_{p-1}^{-}( ({\bf z_*}+{\bf z})^{-t}),  \pi_{p+1}^{+}(({\bf z_*}+{\bf z})^t)\right).\cr}$$
Observe that
$${\bf z_*}+{\bf z} \in \Pi_{\pi_{p-1}^{-}( ({\bf z_*}+{\bf z})^{-t}),\pi_{p+1}^{+}(({\bf z_*}+{\bf z})^t),t},$$
which implies that $\Theta$ is injective. Let us prove that it is  $M e^{At}$ bi-Lipschitz if $E_{p-1}^{-}\times E_{p+1}^{+}$ is munished with the supremum norm. The fact that it is $e^{At}$ Lipschitz is an immediate consequence of Lemma 2. To prove that the inverse is $M e^{At}$ Lipschitz, one can note that for every ${\bf z}$, ${\bf z}'$ in $E^i{}^{\perp}$, one has $\left({\bf z}^*+{\bf z}, {\bf z}^*+{\bf z}'\right)\not\in{\cal W}$. So
either
$$\left({\bf z}^*+{\bf z}, {\bf z}^*+{\bf z}'\right)\in \mathrm{Cl}({\cal W}_{p-1}^{-})$$ or
$$\left({\bf z}^*+{\bf z}, {\bf z}^*+{\bf z}'\right)\in \mathrm{Cl}({\cal W}_{p+1}^{+}).$$
One deduces that
$$\left(({\bf z}^*+{\bf z})^{-t}, ({\bf z}^*+{\bf z}')^{-t}\right)\in {\cal W}_{p-1}^{-}$$ or
$$\left(({\bf z}^*+{\bf z})^t, ({\bf z}^*+{\bf z}')^t\right)\in {\cal W}_{p+1}^{+},$$ which implies that
$$\Vert \pi_{p-1}^{-}(({\bf z}^*+{\bf z})^{-t})- \pi_{p-1}^{-}(({\bf z}^*+{\bf z}')^{-t})\Vert\geq M^{-1} e^{-At} \Vert {\bf z}-{\bf z}'\Vert$$ or
$$\Vert \pi_{p+1}^{+}(({\bf z}^*+{\bf z})^{-t})- \pi_{p+1}^{+}(({\bf z}^*+{\bf z}')^{-t})\Vert\geq M^{-1} e^{-At} \Vert {\bf z}-{\bf z}'\Vert$$ 
Consequently, one knows that $\Theta$ is a homeomorphism and that $\Pi_{{\bf z}, {\bf z}',t}$ is the graph of a function $\theta_{i,t}: E^i\to E^{i}{}^{\perp}$. By Proposition 4, the fact that
$\left((\Pi_{{\bf z}, {\bf z}',t})^s\times (\Pi_{{\bf z}, {\bf z}',t})^s\right)\setminus \mathrm{diag}\subset W_p$,  for every $s\in[-t,t]$, implies that 
this graph is $N_t$ Lipschitz . We can define a similar map $\theta'_{i,t}: E'{}{^i}\to E'{}^{i}{}^{\perp}$ and it will also be $N_t$ Lipschitz. \hfill$\Box$

\bigskip
For every $t\geq 0$, define 
$$\Pi_{p,t}=\Gamma_{\psi_{{\bf 0}}^{t}}\cap \Gamma'_{\psi'{}_{{\bf 0}}^{-t}}= \left(E_{p}^{+}\right)^t\cap \left(E_{p}^{-}\right)^{-t}.$$

Each graph $\Gamma_{\psi_{{\bf 0}}^{t}}$ and $\Gamma'_{\psi_{{\bf 0}}^{-t}}$ contains $\{{\bf 0}\}\cup\Sigma_p$ because this is the case for $E_{p}^{+}$ and $E_{p}^{-}$: one deduces that $\Pi_{p,t}$ contains  $\{{\bf 0}\}\cup\Sigma_p$.

Each graph $\Gamma_{\psi_{{\bf 0}}^{t}}$ and $\Gamma'_{\psi_{{\bf 0}}^{-t}}$ is invariant by $\varphi$ because this is the case for $E_{p}^{+}$ and $E_{p}^{-}$ and because $\xi$ is invariant by $\varphi$: one deduces that $\Pi_{p,t}$ is invariant by $\varphi$.

For every $t\geq 1$ and every $i\in\Z$ , one may write $\Pi_{p,t}$ as the graph of a function $\theta_{i,t}:E^i\to E^i{}^{\perp}$ and as the graph of function $\theta_{i,t}:E'{}^i\to E'{}^i{}^{\perp}$ and all these functions are $N_t$ Lipschitz.

\medskip
Using Ascoli's theorem, one may find  for ever $i\in\Z$ 
 a map $\theta_i:E^i\to E^i{}^{\perp}$, a map $\theta'_i:E'{}^i\to E'{}^i{}^{\perp}$, both $N_t$ Lipschitz for every $t>0$, and a sequence $(t_n)_{n\geq 0}$ satisfying $\lim_{n\to+\infty}t_n=+\infty$, such that each sequence $(\theta_{i,t_n})_{n\geq 0}$ converges to $\theta_i$ and each sequence $(\theta'_{i,t_n})_{n\geq 0}$ converges to $\theta'_i$ for the compact-open topology. The graphs of $\theta_i$ and $\theta'_i$ are equal and independant of $i$. The topological plane $\Pi_p$ obtained in that way is invariant by $\varphi$, contains  $\{{\bf 0}\}\cup\Gamma_p$ and satisfies 
$$\left((\Pi_p)^t\times (\Pi_p^t)\right)\setminus\mathrm{diag}\subset{\cal W}_p$$ for every $t\in\R$. \hfill$\Box$

\bigskip
\noindent{\it Proof of Proposition 2.} \enskip  Let $\Pi_p$ be a plane given by Proposition 5. The curve $\Sigma_p$ bounds a disk $\Delta_{p}\subset \Pi_p$. This  disk is invariant by $\varphi$ because it is the case for $\Pi_p$ and $\Sigma_p$, so the assertion {\bf i)} is satisfied. The assertions {\bf iii)} and {\bf iv)} being true on $\Pi_p$ are of course true on $\Delta_p$. The map $Q_1$ sends  homeomorphically $\Pi_p$ onto $\R^2$ and satisfies $Q_1({\bf 0})=0$. Moreover its sends $\Sigma_p$ onto $S_{p/q}$, which implies that it sends $\Delta_{p}$ onto $D_{p/q}$. Consequently one has ${\bf 0}\in \Delta_p$, which means that {\bf i)} is true.

Let us prove now that for every ${\bf z}\in \Delta_{p}\setminus (\{{\bf 0}\}\cup\Sigma_p)$, one has 
$$\lim_{t\to-\infty} {\bf z}^t={\bf 0}, \enskip \lim_{t\to+\infty} d({\bf z}^t, \Sigma_p)=0.$$Observe that $\Sigma_p$ bounds the disk $(\Delta_{p})^t\subset (\Pi_p)^t$, which is also sent  homeomorphically on $D_{p/q}$ by each $Q_i$ and $Q'_i$. Consequently, the orbit of every point ${\bf z}\in \Delta_{p}\setminus (\{0\}\cup\Sigma_q)$ is bounded. The flow $\xi$ being a gradient flow, the sets $\alpha({\bf z})$ and $\omega({\bf z})$ are not empty, either reduced to ${\bf 0}$ or included in one of the curves of singularities. The only possible circle is $\Sigma_p$, because $({\bf 0},{\bf z}^t)\in W_p$ for every $t\in\R$ and $({\bf 0},{\bf z}')\in W_{p'}$ for every ${\bf z}'\in\Sigma_{p'}$. Using the fact that ${\bf h}( {\bf z}')>{\bf h}({\bf 0})$ if ${\bf z}'\in \Sigma_p$, one deduces that {\bf vi)} is satisfied.

It remains to prove {\bf v)}, which means that $\Delta_{p}$ is invariant by the flow. It is sufficient to prove that the disk $\Delta_p$ is independent of the plane $\Pi_p$  given by Proposition 5. Indeed $(\Pi_p)^t$ also satisfies the properties stated in Proposition 5, and $\Sigma_p$ bounds the disk $(\Delta_{p})^t\subset (\Pi_p)^t$. We will give a proof by contradiction, supposing that $\Pi_p$ and $\Pi'_p$ are two planes given by Proposition 5, such that $\Sigma_p$ bounds two distinct disks $\Delta_p\subset \Pi_p$ and $\Delta'_p\subset \Pi'_p$. 

Suppose that ${\bf z}\in \Delta_{p}\setminus\Delta'_{p}$. The map $Q_1$ inducing homeomorphisms from $\Delta_{p}$ and $\Delta'_{p}$  onto $D_{p/q}$, there exists ${\bf z}'\in \Delta'_{p}$ such that $Q_1({\bf z}')=Q_1({\bf z})$, which implies that $({\bf z}, {\bf z}')\not\in {\cal W}$. Therefore, one has $({\bf z}, {\bf z}')\in \mathrm{Cl}({\cal W}_{p+1}^{+})$ or $({\bf z}, {\bf z}')\in \mathrm{Cl}({\cal W}_{p-1}^{-})$. In the first case, there exist $p'>p$ and $T>0$ such that $({\bf z}^{t}, {\bf z}'{}^{t})\in {\cal W}_{p'}$ for every $t>T$. In the second case, there exist $p'<p$ and $T<0$ such that $({\bf z}^{t}, {\bf z}'{}^{t})\in {\cal W}_{p'}$ for every $t<T$. Let us consider the second case. The fact that $\lim_{t\to -\infty} {\bf z}^t ={\bf 0}$ and that $({\bf 0}, {\bf z}^t)\in {\cal W}_p$ implies that the line generated by ${\bf z}^t$ is approaching the space $E_{p}({\bf 0})$ when $t$ tends to $-\infty$. We have a similar result for   ${\bf z}'{}^t$. The fact that $({\bf z}^{t+s}, {\bf z}'{}^{t+s})\in {\cal W}_{p'}$ for every $s\in(-\infty, -t+T)$ implies that the line generated by ${\bf z}^{t} -{\bf z}'{}^{t}$ is approaching the space $E_{p'}({\bf 0})$ when $t$ tends to $-\infty$. 

As we know that the eigenvalues of $D\xi({\bf 0})\vert_{E_{p'}({\bf 0})}$ are smaller than 
the eigenvalues of $D\xi({\bf 0})\vert_{E_{p}({\bf 0})}$, we deduce that there exist $C>0$, $C'>0$ and $\mu>\mu'$  such that, for $-t$ large enough, one has  
$$\Vert {\bf z}^{t}\Vert \leq C e^{\mu t}, \enskip  \Vert {\bf z}'{}^{t}\Vert \leq C e^{\mu t},Ê\enskip \Vert{\bf z}^{t} -{\bf z}'{}^{t}\Vert\geq C' e^{\mu' t},$$ which is impossible.

Let us consider the first case. The fact that $\Sigma_p\times\Sigma_p\setminus\mathrm{diag}$ is included in ${\cal W}_p$ implies that $0$ is an eigenvalue of $D\xi({\bf z}'')\vert_{E_{p}({\bf z}'')}$ for every ${\bf z}''\in \Sigma_q$. Consequently, there exists a uniform upper bound $\mu'>0$ of the spectrum of $D\xi({\bf z}'')\vert_{E_{p'}({\bf z}'')}$. The fact that $\omega({\bf z})$ and $\omega({\bf z'})$ are included in $\Sigma_p$ and that $({\bf z}^{t}, {\bf z}'{}^{t})\in {\cal W}_{p'}$ for every $t>T$ implies that $\lim_{t\to +\infty} \Vert{\bf z}^{t} -{\bf z}'{}^{t}\Vert =0$.  Moreover if $t$ is large, then for every $s\in(-T-t,+\infty)$ one has $({\bf z}^{t+s}, {\bf z}'{}^{t+s})\in {\cal W}_{p'}$. This implies that if $t$ is large and ${\bf z}'{}^{t}$ is close to a point ${\bf z}''\in\Sigma_p$, then it is also the case for ${\bf z}^{t}$ and the line  generated by ${\bf z}^{t} -{\bf z}'{}^{t}$ is close to $E_{p'}({\bf z}'')$. Consequently, if $t$ is large, then 
$$ \Vert{\bf z}^{t+1} -{\bf z}'{}^{t+1}\Vert\geq\Vert{\bf z}^{t} -{\bf z}'{}^{t}\Vert,$$
which is impossible. \hfill$\Box$

\bigskip
\noindent{\bf Remark.} \enskip  The uniqueness property that has been stated in the previous proof permits us to give a construction of $\Delta_p$. Let $\Pi_{p,t}= \left(E_{p}^{+}\right)^t\cap \left(E_{p}^{-}\right)^{-t}$ be the plane defined in the proof of Proposition 5 and $\Delta_{p,t}\subset \Pi_{p,t}$ the disk bounded by $\Sigma_t$, then one has $$\Delta_p=\lim_{t\to+\infty} \Delta_{p,t},$$
for each natural topology (induced by the Hausdorff distance, associated to the $C^0$ topology on maps $\theta_i:E^i\to E^i{}^{\perp}$ or on maps $\theta'_i:E'{}^i\to E'{}^i{}^{\perp}$). The uniqueness property, and consequently its invariance by the flow is a consequence of the fact that $\{0\}\cup S_{p/q}$  is a maximal unlinked fixed points set of $f^q$. The fact that $\xi$ is a gradient flow was not essential in the proof, nevertheless the proof is easier in this case (see Propostion 5.2.1 in \cite{Lec2}).

\bigskip

\noindent{\large {\bf 6. Appendix}}

\bigskip
The goal of the appendix is to give a proof of Proposition 4. This proposition will result from the technical Lemma 6. A particular case of this lemma, that we will explain in details, implies Proposition 3.  The proof of Lemma 6, in this particular case, is nothing but the proof  of Lemma 2.5.1 in \cite{Lec2}. We will look at a wider situation than the one studied in the present paper by looking at a general tridiagonal and monotonically symmetric Lipschitz vector field.

\medskip

We fix an integer $r\geq 2$ and consider the finite dimensional vector space
$$	F= \left\lbrace {\bf x}=(x_{i})_{i\in {\Z}}\in \R^{\Z}\enskip 
\enskip\vert \enskip x_{i+r}=x_{i},\enskip 
\mathrm{ for\enskip all}\enskip i\in {\Z}\right\rbrace.
$$
Unlike in Section 2, we munish $F$ with the sup norm $\Vert\enskip\Vert$ where $\Vert {\bf x}\Vert =\max_{i\in\Z} \vert x_i\vert$. 

\medskip
We consider on $F$ a tridiagonal and monotonically symmetric Lipschitz vector field $\zeta=(\zeta_i)_{i\in\Z}$. More precisely we suppose that
 $$
\zeta
_{i}({\bf x})=\zeta_i(x_{i-1}, x_i,x_{i+1}),$$where

\smallskip
\noindent-\enskip\enskip the map $x\mapsto \zeta_i(x_{i-1}, x,x_{i+1})$ is $K$ Lipschitz;

\smallskip
\noindent-\enskip\enskip the maps  $x\mapsto \zeta_i(x, x_i,x_{i+1})$ and $x\mapsto \zeta_i(x_{i-1}, x_i,x)$ are $K$ bi-Lipschitz homeomorphisms;

\smallskip
\noindent-\enskip\enskip if $x\mapsto \zeta_i(x_{i-1}, x_i,x)$ is increasing, then  $x\mapsto \zeta_{i+1}(x, x_{i+1},x_{i+2})$ is increasing and we set $\sigma_i=1$;

\smallskip
\noindent-\enskip\enskip if $x\mapsto \zeta_i(x_{i-1}, x_i,x)$ is decreasing, then  $x\mapsto \zeta_{i+1}(x, x_{i+1},x_{i+2})$ is decreasing and we set $\sigma_i=-1$.

\smallskip

Here again, we write ${\bf x}^t$ for the image at time $t$ of a point ${\bf x}$ by the flow.

\medskip
We will prove the following:

\bigskip
\noindent{\sc Lemma 6:}\enskip\enskip{\it For every integer $l\geq 2$ and  every real number $\mu\in(0,1]$, there exists:

\smallskip
\noindent-\enskip\enskip  a sequence of polynomials $(P_k)_{k\geq 0} $ in $\R[X_1,X_2]$ depending on $K$, $l$ and $\mu$, such that
$P_k(0,X_2) = a_kX_2^k$ with $a_k>0$ if $k\leq [l/2]+1$;

\smallskip
\noindent-\enskip\enskip  a sequence of polynomials $(Q_k)_{k\geq 0} $ in $\R[X_1,X_2]$ depending on $K$, $l$, and $\mu$, such that
$Q_k(0,X_2) = b_kX_2^k$ with $b_k>0$;

\smallskip 
\noindent-\enskip\enskip 
real numbers $\varepsilon_0>0$ and $t_0>0$, 

\smallskip
\noindent which verify the following:

\medskip

\noindent for every ${\bf x}$, ${\bf x}'$ such that
$$\eqalign{ \vert x_i-x'_i\vert  &\geq \mu\Vert {\bf x}-{\bf x}'\Vert \cr
\vert x_{i+1}-x'_{i+1}\vert &\leq \varepsilon\Vert {\bf x}-{\bf x}'\Vert\cr
\dots&\dots\dots\cr
\vert x_{i+l}-x'_{i+l}\vert &\leq \varepsilon\Vert {\bf x}-{\bf x}'\Vert\cr 
\vert x_{i+l+1}-x'_{i+l+1}\vert  &\geq\mu\Vert {\bf x}-{\bf x}'\Vert\cr}$$
with $\varepsilon \leq \varepsilon_0$, and for every $t\in[-t_0,t_0]\setminus\{0\}$, one has
$$\eqalign{\sigma'_{i+k}\, \,\mathrm{sign}(t)^k(x_{i+k}^t-x'{}^t_{i+k})&\geq P_k(\varepsilon, \vert t\vert)\,\Vert {\bf x}-{\bf x}'\Vert\cr
\vert x_{i+k+1}^t-x'{}^t_{i+k+1}\vert &\leq Q_k(\varepsilon, \vert t\vert )\,\Vert {\bf x}-{\bf x}'\Vert\cr
\dots&\dots\dots\cr
\vert x_{i+l-k}^t-x'{}^t_{i+l-k}\vert &\leq Q_k(\varepsilon, \vert t\vert)\,\Vert {\bf x}-{\bf x}'\Vert\cr 
\sigma'_{i+l+1-k}\,\,\mathrm{sign}(t)^k(x_{i+l+1-k}^t-x'{}^t_{i+l+1-k})&\geq P_k(\varepsilon, \vert t\vert)\,\Vert {\bf x}-{\bf x}'\Vert,\cr}$$
where
$$\sigma'_{i+k}= \mathrm{sign}(x_i-x'_i)\,\sigma_i\dots \sigma_{i+k-1}$$ and
$$\sigma'_{i+l+1-k}= \mathrm{sign}(x_{i+l+1}-x'_{i+l+1})\,\sigma_{i+l+1-k}\dots \sigma_{i+l}.$$

\medskip
\noindent Moreover if $l=2m-1$ is odd and $\sigma'_{i+m-1}\sigma_{i+m-1}= \sigma'_{i+m+1}\sigma_{i+m}$, then
$$\sigma'_{i+m}\,\,\mathrm{sign}(t)^m( x_{i+m}^t-x'{}^t_{i+m})\geq 2P_{m}(\varepsilon, \vert t\vert)\,\Vert {\bf x}-{\bf x}'\Vert,$$}
where 
$$\sigma'_{i+m}= \sigma'_{i+m-1}\sigma_{i+m-1}= \sigma'_{i+m+1}\sigma_{i+m}.$$

\bigskip
 
\noindent{\it Proof.}\enskip\enskip By replacing the vector field $\xi$ with $-\xi$, one changes every $\sigma_i$ into $-\sigma_i$. Thus it is sufficient to prove the lemma for $t>0$.
Let us begin by defining two sequences $(a_k)_{k\geq 0}$ and $(b_k)_{k\geq 0}$ by the induction equations
$$ a_{k+1}= k^{-1}\left(K^{-1} a_k- 2K b_k\right), \,b_{k+1}= 3Kk^{-1} b_k$$
and the initial conditions
$$a_0= \mu/2,\enskip b_0=\delta,$$
where $\delta>0$ is chosen sufficiently small to ensure that $a_k$ is positive if $k\leq [l/2]+1$. The number $\delta$ depends on $K$, $l$ and $\mu$.  Suppose moreover than $\delta\leq \mu$ and  write
$$\varepsilon_0={\delta\over 2}, \, t_0={\log (\delta/2)\over 3K}.$$

Let us continue by defining two sequences of polynomials $(P_k)_{k\geq 0}$, $(Q_k)_{k\geq 0}$ in $\R[X_1,X_2]$ by the induction equations
$${\partial\over \partial X_2} P_{k+1} =K^{-1}P_k-2K Q_k,\enskip P_{k+1}(X_1,0)= -X_1,$$
$${\partial\over \partial X_2} Q_{k+1}  =3KQ_k(X_1,X_2),\enskip Q_{k+1}(X_1,0)= X_1,$$
and the initial conditions
$$P_0=\mu/2,  \enskip Q_0=\delta.$$
Observe that
$$P_k(0,X_2)= a_kX_2^k, \enskip Q_k(0,X_2)= b_kX_2^k.$$
Note also that $Q_k(x_1,x_2)\geq 0$ if $x_1\geq 0$ and $x_2\geq 0$.

\medskip
We will prove by induction on $k\in\{0,\dots, [l/2]\}$ that for every $t\in(0,t_0]$, one has
$$\eqalign {\sigma'_{i+k}( x_{i+k}^t-x'{}^t_{i+k})&\geq P_k(\varepsilon, t)\,\Vert {\bf x}-{\bf x}'\Vert\cr
\sigma'_{i+l+1-k}( x_{i+l+1-k}^t-x'{}^t_{i+l+1-k})&\geq P_k(\varepsilon, t)\,\Vert {\bf x}-{\bf x}'\Vert\cr}$$ and 
$$\vert x_{j}^t-x'{}^t_{j}\vert\leq Q_k(\varepsilon, t)\,\Vert {\bf x}-{\bf x}'\Vert$$
if $i+k<j< i+l+1-k$.

\medskip
By hypothesis, the vector field is $3K$ Lipschitz. Using Gronwall's Lemma, like in Lemma 2, one deduces that
$$e^{-3K t}\Vert{\bf x}^t-{\bf x}'{}^t\Vert \leq\Vert {\bf x}^t-{\bf x}'{}^t\Vert \leq  e^{3Kt}\Vert{\bf x}-{\bf x}'\Vert .$$ 
Therefore, for every $i\in\Z$, one has
$$\vert x_i -x'_i\vert \leq e^{3Kt}\Vert{\bf x}-{\bf x}'\Vert ,$$
which implies
$$\vert \dot x_i -\dot x'_i\vert \leq 3Ke^{3Kt}\Vert{\bf x}-{\bf x}'\Vert ,$$
and so
$$\vert (x_i^t -x'{}_i^t)-(x_i-x'_i)\vert \leq e^{3Kt}\Vert{\bf x}-{\bf x}'\Vert .$$

By definition of $\varepsilon_0$ and $t_0$ one deduces that the induction hypothesis is satisfied for $k=0$.

\medskip

Suppose that the induction hypothesis has been proven until $k$, where $k<[l/2]$ and let us prove it for $k+1$.
One has $$\eqalign{\dot x_{i+k+1} -\dot x'_{i+k+1}
&=\enskip \enskip \zeta_{i+k+1}( x_{i+k}, x_{i+k+1}, x_{i+k+2}) -\zeta_{i+k+1}( x'_{i+k}, x_{i+k+1}, x_{i+k+2}) \cr
&\enskip\enskip +\zeta_{i+k+1}( x'_{i+k}, x_{i+k+1}, x_{i+k+2}) -\zeta_{i+k+1}( x'_{i+k}, x'_{i+k+1}, x_{i+k+2}) \cr
&\enskip\enskip+\zeta_{i+k+1}( x'_{i+k}, x'_{i+k+1}, x_{i+k+2}) -\zeta_{i+k+1}( x'_{i+k}, x'_{i+k+1}, x'_{i+k+2}) \cr}.$$
By definition of  $\sigma_{i+k+1}$, $\sigma'_{i+k+1}$ and by hypothesis, one deduces that
$$\sigma_{i+k+1}'(\dot x_{i+k+1}^t-\dot x'{}^t_{i+k+1})\geq \left(K^{-1} P_k(\varepsilon, t) -2K Q_k(\varepsilon, t)\right) \Vert{\bf x}-{\bf x}'\Vert.$$
Consequently, one has 
$$\sigma'_{i+k+1}( (x_{i+k+1}^t -x'{}^t_{i+k+1})-(x_{i+k+1}-\dot x'_{i+k+1}))\geq\left(\int_0^t (K^{-1} P_k(\varepsilon, u)- 2 KQ_k(\varepsilon, u))\,du\right) \Vert{\bf x}-{\bf x}'\Vert ,$$
which implies that
$$\eqalign{\sigma'_{i+k+1}(x_{i+k+1}^t-x'{}^t_{i+k+1})&\geq \left(\int_0^t K^{-1} P_k(\varepsilon, u) -2K Q_k(\varepsilon, u)\, du-\varepsilon\right )\Vert{\bf x}-{\bf x}'\Vert\cr&=
 P_{k+1}(\varepsilon, t)\Vert{\bf x}-{\bf x}'\Vert.\cr}$$
 The same proof gives us
 $$\sigma_{i+l-k}(x_{i+l-k}^t-x'{}^t_{i+l-k})\geq
 P_{k+1}(\varepsilon, t)\,\Vert{\bf x}-{\bf x}'\Vert.$$

\noindent Now, let us fix $j\in\{k+2,\dots, l-k-1\}$. One has

$$\vert\dot x_{j}^t-\dot x'{}^t_{j}\vert\leq 3K Q_k(\varepsilon, t) \Vert{\bf x}-{\bf x}'\Vert$$
and so
$$\vert(x_{j}^t-x'{}_{j}^t)-(x_i-x'_i))\vert\leq \left(\int_0^t 3K Q_k(\varepsilon, u)\, du\right)\Vert{\bf x}-{\bf x}'\Vert$$
which implies
$$\vert x_{j}^t-x'{}^t_{j}\vert \leq \left(\int_0^t 3K Q_k(\varepsilon, u)\, du+\varepsilon\right )\Vert{\bf x}-{\bf x}'\Vert=Q_{k+1}(\varepsilon, t) \Vert{\bf x}-{\bf x}'\Vert$$

It remains to study the case where $l=2m-1$ is odd and $\sigma'_{i+m-1}\sigma_{i+m-1}= \sigma'_{i+m+1}\sigma_{i+m}$. One has
$$\sigma'_{i+m}(\dot x_{i+m}^t-\dot x'{}^t_{i+m})\geq \left(2K^{-1} P_{m-1}(\varepsilon, t) -K Q_{m-1}(\varepsilon, t)\right) \Vert{\bf x}-{\bf x}'\Vert$$
which implies
$$\eqalign{\sigma'_{i+m}(x_{i+m}^t-x'{}^t_{i+m})&\geq \left(\int_0^t 2K^{-1} P_{m-1}(\varepsilon, u) -K Q_{m-1}(\varepsilon, u)\, du-\varepsilon\right )\Vert{\bf x}-{\bf x}'\Vert\cr
&\geq \left(\int_0^t 2K^{-1} P_{m-1}(\varepsilon, u) -4K Q_{m-1}(\varepsilon, u)\, du-2\varepsilon\right )\Vert{\bf x}-{\bf x}'\Vert\cr
&= 2P_{m}(\varepsilon, t)\Vert{\bf x}-{\bf x}'\Vert\cr}$$
\hfill$\Box$

\bigskip
Let us explain why  Proposition 3 can be deduced from Lemma 6 in the case where $\varepsilon=0$. The function defined by the formula 
$${\cal L}({\bf x})=\sum_{0<i\leq r} \sigma_i \,\mathrm{sign}( x_i) \,\mathrm{sign}( x_{i+1})$$
extends naturally to the open set 
$$U=\{{\bf x}\in F\enskip\vert\enskip x_i=0\Rightarrow \sigma_{i-1}\sigma_i \,x_{i-1}x_{i+1}<0\}.$$
Supposing $\varepsilon=0$, Lemma 6 permits to determine the sign of every $x^t_j-x'{}^t_j$, $j\in\{i, \dots, i+l+1\}$, $\vert t\vert\in(0,t_0]$,  with the exception of $x^t_{i+m}-x'{}^t_{i+m} $, in the case where $j=2m-1$ is odd and $\sigma'_{i+m-1}\sigma_{i+m-1}\not= \sigma'_{i+m+1}\sigma_{i+m}$. One has
$$\eqalign{\mathrm{sign } (x^t_{i+k}-x'{}^t_{i+k})&=\mathrm{sign} (t)^k\,\sigma'_{i+k}\cr
\mathrm{sign } (x^t_{i+l+1-k}-x'{}^t_{i+l+1-k})&=\mathrm{sign} (t)^k\,\sigma'_{i+l+1-k}\cr}.$$
Recall that $\sigma'_{i+k}=\sigma_{i+k-1}\sigma'_{i+k-1}$ and $\sigma'_{i+l+1-k}=\sigma_{i+l+1-k}\sigma'_{i+l+2-k}$, if $1\leq k\leq [l/2]$, which implies that 
$$\mathrm{sign}( x_j^t-x'{}^t_j) \,\mathrm{sign}( x_{j+1}^t-x'{}^t_{j+1})= \mathrm{sign} (t)\,\sigma_j \,\sigma'_j\sigma'_{j+1} =\mathrm{sign} (t),$$ if $j\in\{i,\dots, i+[l/2]-1, i+[l/2]+1, \dots, i+l\}$. 

\medskip
In the case where $t$ is positive, the value of
$$\sum_{i\leq j\leq i+l} \sigma_j \,\mathrm{sign}( x_j^t-x'{}^t_j) \,\mathrm{sign}( x_{j+1}^t-x'{}^t_{j+1})$$
 is
$$\cases {l+\sigma_{i+m} \,\sigma'_{i+m}\sigma'_{i+m+1}& if $l=2m$ is even,\cr
l+1& if $l=2m-1$ is odd and $\sigma'_{i+m-1}\sigma_{i+m-1}= \sigma'_{i+m+1}\sigma_{i+m}$,\cr
l-1& if $l=2m-1$ is odd and $\sigma'_{i+m-1}\sigma_{i+m-1}\not= \sigma'_{i+m+1}\sigma_{i+m}$.\cr}$$
In the case where $t$ is negative, it is
$$\cases {-l+\sigma_{i+m} \,\sigma'_{i+m}\sigma'_{i+m+1}& if $l=2m$ is even,\cr
-l-1& if $l=2m-1$ is odd and $\sigma'_{i+m-1}\sigma_{i+m-1}= \sigma'_{i+m+1}\sigma_{i+m}$,\cr
-l+1& if $l=2m-1$ is odd and $\sigma'_{i+m-1}\sigma_{i+m-1}\not= \sigma'_{i+m+1}\sigma_{i+m}$.\cr}$$
 The difference between the two values is a positive multiple of $4$, except in the case where $l=1$ and $\sigma'_{i}\sigma_{i}\not= \sigma'_{i+2}\sigma_{i+1}$, which means that $\sigma_i \mathrm{sign}(x_i-x'_i)\not= \sigma_{i+1} \mathrm{sign}(x_{i+2}-x'_{i+2})$. 

Consequently, if
${\bf x}$, 
${\bf x}'$ are two distinct points of $F$ such that ${\bf x}'-{\bf x}\not\in U$, there exists $t_0 >0$ such
that
${\bf x'}^{t}-{\bf x}^{t}\in U$ if  $0< \left \vert t\right 
\vert \leq t_0$ and one has
$${\cal L}({\bf x'}^{t_0} -{\bf x}^{t_0
})>{\cal L}({\bf x'}^{-t_0 } -{\bf x}^{-t_0 })\enskip.
$$
The proof given above in the case where $\varepsilon=0$ is nothing but the proof of Lemma 2.5.1 in \cite{Lec2}.
It implies Proposition 3 when applied to the vector field $\xi$ because $L={\cal L}/4$. 

\medskip
Lemma 6 tells us more. For every $l\geq 2$  and $\mu$, there exists $t_0$ such that for every $t\in(0,t_0]$, there exists $\epsilon$ such that $P_k(\varepsilon, t)> 0$ for every $k\leq [l/2]+1$. Indeed one has  $P_k(0, t)= a_kt^k> 0$. Consequently, if ${\bf x}$ and ${\bf x'}$ satisfies 
$$\eqalign{ \vert x_i-x'_i\vert  &\geq \mu\Vert {\bf x}-{\bf x}'\Vert \cr
\vert x_{i+1}-x'_{i+1}\vert &\leq \varepsilon\Vert {\bf x}-{\bf x}'\Vert\cr
\dots&\dots\dots\cr
\vert x_{i+l}-x'_{i+l}\vert &\leq \varepsilon\Vert {\bf x}-{\bf x}'\Vert\cr 
\vert x_{i+l+1}-x'_{i+l+1}\vert &\geq\mu\Vert {\bf x}-{\bf x}'\Vert\cr}$$
one can determine the signs of  $x^t_j-x'{}^t_j$, $j\in\{i, \dots, i+l+1\}$, with the exception of $x^t_{i+m}-x'{}^t_{i+m} $, in the case where $j=2m-1$ is odd and $\sigma'_{i+m-1}\sigma_{i+m-1}\not= \sigma'_{i+m+1}\sigma_{i+m}$.  Like in the case where $\varepsilon =0$, 
the value of
$$\sum_{i\leq j\leq i+l} \sigma_j \,\mathrm{sign}( x_j^t-x'{}^t_j) \,\mathrm{sign}( x_{j+1}^t-x'{}^t_{j+1})$$
is larger at time $t$ than at time $-t$ except if $l=1$ and $\sigma_i \mathrm{sign}(x_i-x'_i)\not= \sigma_{i+1} \mathrm{sign}(x_{i+2}-x'_{i+2})$.

\medskip
As a consequence, if one defines, for $t>0$, the sets
$${\cal X}_t =\{({\bf x}, {\bf x}')\in F^2\, \vert  \enskip \vert s\vert\leq t\Rightarrow {\bf x}^s-{\bf x}'{}^s\in U \}$$
and $${\cal Y}_t =\left\{{{\bf x}-{\bf x}'\over \Vert {\bf x}-{\bf x}' \Vert }\, \vert  \enskip ({\bf x}, {\bf x}')\in{\cal X}_t \right\},$$
one deduces that the closure of ${\cal Y}_t$ is included in $U$. This implies Proposition 4.

\bigskip
Let us conclude this appendix with the following result:

\bigskip
\noindent{\sc Lemma 7:}\enskip\enskip{\it The space $\mathrm{Diff}^1_{**}({\bf D})$ of area preserving $C^1$ diffeomorphisms of $\bf D$ that fix $0$ and every point of $\bf S$ is path-connected when munished with the $C^1$ topology.}

\bigskip
\noindent{\it Proof.}\enskip\enskip Let us begin by proving that the space $\mathrm{Diff}^1_{*}({\bf D})$ of area preserving $C^1$ diffeomorphisms of $\bf D$ that fix every point of $\bf S$ is path-connected when munished with the $C^1$ topology. Let us consider the symplectic polar system of coordinates $(\theta, \lambda)\in\T^1\times\R$ defined on on $\R^2\setminus\{0\}$ by
$$\eqalign{\lambda&=x^2+y^2\cr (x,y)&=(\sqrt \lambda \cos 2\pi\theta, \sqrt \lambda \sin 2\pi\theta)\cr}$$ Every $\Phi\in\mathrm{Diff}^1_{*}({\bf D}) $ induces a $C^1$ area preserving diffeomorphism, defined in a neighborhood of $\T^1\times\{1\}$ in $\T^1\times(0,1]$ that fixes every point of $\T^1\times\{1\}$. The image of a circle $\lambda =\lambda_0$ is a graph $\lambda =\psi( \theta)$ if $\lambda_0$ is close to $1$. This implies that $\Phi$ is defined by a generating function in a neighborhood of $\T^1\times\{1\}$: if one writes $\Phi(\theta,\lambda)=(\Theta, \Lambda)$, one knows that $(\Theta, \lambda)$ defines a system of coordinates  in a neighborhood of $\T^1\times\{1\}$ in $\T^1\times(0,1]$  and that $(\theta-\Theta) d\lambda+(\Lambda-\lambda) d\Theta$ is a $C^1$ exact form. There exists a  $C^2$ function $S $ defined on an annulus $\T^1\times[1-\eta_0,1\}$ such that $$\theta-\Theta={\partial S\over\partial \lambda}, \enskip \Lambda-\lambda ={\partial S\over\partial \Theta}$$
and such that $${\partial ^2 S\over\partial\Theta\partial \lambda}>-1.$$
The circle $\lambda=1$ being invariant and contained in the fixed point set of $\Phi$, one has
$${\partial S\over\partial \Theta}(\Theta, 1)= {\partial S\over\partial \lambda}(\Theta, 1)=0$$ which implies
$${\partial ^2 S\over\partial\Theta\partial \lambda}(\Theta, 1)=0.$$

Let $\nu:[0,+\infty)\to[0,1]$ be a $C^2$ function such that 
$$\nu(t)=\cases {1 & if $t\in[0,1/3]$,\cr
0 & if $t\geq 2/3$.\cr}$$
For every $\eta\in(
0,\eta_0]$ define
$$S_{\eta}: (\Theta, \lambda)\mapsto\nu\left({1-\lambda\over \eta}\right) S(\Theta, \lambda).$$
 
 One has
 $$ {\partial ^2 S_{\eta}\over\partial\Theta\partial \lambda}= \nu\left({1-\lambda\over \eta}\right) {\partial ^2 S\over\partial\Theta\partial \lambda} -{1\over \eta} \nu'\left({1-\lambda\over \eta}\right){\partial S\over\partial\Theta} .
 $$
 The fact that $\displaystyle{\partial ^2 S\over\partial\Theta\partial \lambda}(\Theta, 1)=0$ implies that the quantity $\displaystyle{1\over \eta}{\partial S\over\partial\Theta}$ tends to zero uniformly on the annulus of equation $1-\eta \leq \lambda\leq 1$, when $\eta$ tends to $0$. As a consequence, if $\eta$ is sufficiently small, one has $${\partial ^2 S\over\partial\Theta\partial \lambda}>-1$$ on $\T^1\times [1-\eta_0, 1]$. The equations
 $$\theta-\Theta=s{\partial S_{\eta}\over\partial \lambda}, \enskip \Lambda-\lambda =s{\partial S_{\eta}\over\partial \Theta}$$ defines a family $(\Phi_{\eta,s})_{s\in[0,1]}$ of area preserving $C^1$ diffeomorphisms on the annulus of equation $1-\eta_0\leq \lambda\leq 1$ that fix all point in a neighborhood of the circle $\lambda =1-\eta_0$. Extending these maps by the identity on ${\bf D}$, one gets a continous arc in $\mathrm{Diff}^1_*({\bf D})$ that joins the identity to a diffeomorphism $\Phi_{\eta,1}$ that coincides with $\Phi$ in a neighborhood of $\bf S$. One can write $\Phi= \Phi'\circ \Phi_{\eta,1}$, where $\Phi'$ coincides with the identity in a neighborhood of $\bf S$. It remains to prove that $\Phi'$ can be joined to the identity in $\mathrm{Diff}^1_*({\bf D})$. A classical result of Zehnder \cite{Z}, whose proof uses generating functions states, that every $C^1$ symplectic diffeomorphism  on a compact symplectic manifold is the limit, for the $C^1$ topology of a sequence of $C^{\infty}$ symplectic diffeomorphisms. If $\cal U$ is a neighborhood of the identity in $\mathrm{Diff}^1_*({\bf D})$, Zehnder's proof permits us to write $\Phi'=\Phi''\circ \Phi'''$, where $\Phi''$ belongs to $\cal U$ and $\Phi'''$ is an area preserving $C^{\infty}$ diffeomorphism of $\bf D$, both of them coinciding with the identity in a neighborhood of $\bf S$.  If $\cal U$ is sufficiently small, $\Phi''$ will be defined by a generating function:  there exist  a $C^2$ function $h$ defined on the disk $X^2+y^2\leq 1$ and $\eta_1>0$  such that
 $$\Phi''(x,y)=(X,Y)\Leftrightarrow\cases{ x=\partial h/\partial y (X,y)\cr Y=\partial h/\partial X (X,y)\cr},$$ where
  $${\partial ^2 h\over\partial X\partial y}>0$$ and
  $$ h(X,y)=Xy \enskip \mathrm{if} \enskip 1-\eta_1 \leq X^2+y^2\leq 1.$$

  Writing $h_s(X,y)= sh(X,y)+(1-s)Xy$, for $s\in[0,1]$, the equations  
  $$\Phi_s''(x,y)=(X,Y)\Leftrightarrow\cases{ x=\partial h_s/\partial y (X,y)\cr Y=\partial h_s/\partial X (X,y)\cr}.$$
 define a continous arc in $\mathrm{Diff}^1_*({\bf D})$ that joins the identity to $\Phi''$. The fact that $\Phi'''$ can be joined to the identity comes from the following classical result whose proof uses a similar result, due to Smale, in the non area preserving case and the classical Moser's Lemma on volume forms (see \cite {Ba} for example): the space of area preserving $C^{\infty}$ diffeomorphisms of $\bf D$, that coincides with the identity on a neighborhood of $\bf S$ is path-connected when munished with the $C^{\infty}$ topology.

 \medskip
 
 Let us prove now that $\mathrm{Diff}^1_{**}({\bf D})$ is path-connected. It is sufficient to prove that for every compact set $\Xi$ in the interior of ${\bf D}$, there exists a continuous family $(\Psi^z)_{z\in\Xi}$ in $\mathrm{Diff}^1_{*}({\bf D})$ such that $\Psi^z(z)=0$. Indeed, every $\Phi\in\mathrm{Diff}^1_{**}({\bf D})$ can be joined to the identity by a path $s\mapsto \Phi_s$ in $\mathrm{Diff}^1_{*}({\bf D})$.  If $\Xi$ is chosen to contain the image of the path $s\mapsto\Phi_s(0)$, then the path $s\mapsto \Psi^{\Phi_s(0)}\circ\Phi_s$ joins $\Phi$ to the identity in $\mathrm{Diff}^1_{**}({\bf D})$.

 For every $r\in(0,1)$ one may find a $C^2$ function $\nu_r:[0,1]\to \R$ such that 
$$\nu(t)=\cases {1 & if $t\in[0,r]$,\cr
0 & if $t\geq 1$.\cr}$$
Denote by $(\Psi_1^s)_{s\in\R}$ and $(\Psi_2^s)_{s\in\R}$ the Hamiltonian flows (for the usual symplectic form $dx\wedge dy$) induced by the function $$H_1:(x,y)\mapsto -y\nu_r(x^2+y^2), \enskip H_1:(x,y)\mapsto x\nu_r(x^2+y^2)$$ respectively,  and for every $z=(x,y)\in\R^2$ set $\Psi^z=\Psi_1^x\circ  \Psi_2^y$. Observe that if $x^2+y^2\leq r$,  then $\Psi^z(z)=0$.
 
 \hfill$\Box$

\bibliographystyle{alpha}
\renewcommand{\refname}{

  \end{document}